\documentclass[11pt,twoside,leqno]{aomamlt2e} 
  \pageno{943}
\received{April 13, 2005}
\allowdisplaybreaks

 \renewcommand{\thesection}{\arabic{section}}
\makeatletter
\renewcommand{\@seccntformat}[1]{\csname
the#1\endcsname.\hspace{0.5em}\setcounter{Subsec}{0}\setcounter{Subsubsec}{0}}\makeatother

 \newtheorem{lemma}{Lemma}

\newtheorem{prop}{Proposition}
\theorembodyfont{\upshape}
\newtheorem{defin}{{\it Definition}}
\newcommand\sdemo[1]{\demo{\scshape #1}}
   \begin{document}
\currannalsline{162}{2005} 

 \title{Hypoellipticity and loss of derivatives}

 \acknowledgements{Research was partially supported by
NSF Grant DMS-9801626.}
 \author{J.\ J.\ Kohn}
 
 \institution{Princeton University, Princeton, NJ
\\
\email{kohn@math.princeton.edu}}

 \shorttitle{Hypoellipticity and loss of derivatives}
  
\centerline{(with an Appendix by Makhlouf Derridj and David S. Tartakoff)}
\vskip12pt

\centerline{\it Dedicated to Yum-Tong Siu for his $60^{\rm th}$ birthday.}

\vglue15pt 
\centerline{\bf Abstract}
\vglue12pt

Let $\{X_1,\dots,X_p\}$ be complex-valued vector fields 
in $\mathbb R^n$ and assume that they satisfy the bracket
condition (i.e. that their Lie algebra spans all vector fields). 
Our object is to study the operator $E=\sum X_i^*X_i$, where $X_i^*$ is
the
$L_2$ adjoint of $X_i$. A result of H\"ormander is that when the $X_i$
 are real then $E$ is hypoelliptic and furthemore it is subelliptic (the
restriction of a destribution $u$ to an open set $U$ is ``smoother'' 
then the restriction of $Eu$ to $U$). When the $X_i$ are complex-valued
 if
the bracket condition of order one is satisfied (i.e. if the
 $\{X_i,[X_i,X_j]\}$ span), then we prove that the operator $E$ is still
subelliptic. This is no longer true if brackets of higher order are
needed to span. For each $k\ge1$ we give an example of two complex-valued
vector fields,
$X_1$ and $X_2$, such that the bracket condition of order $k+1$ is 
satisfied and we prove that the operator $E=X_1^*X_1+X_2^*X_2$ is
hypoelliptic but that it is not subelliptic. In fact it ``loses'' $k$
derivatives in the sense that, for each $m$, there exists a distribution
$u$ whose restriction to an open set $U$ has the property that the
$D^\alpha Eu$ are bounded on $U$ whenever $|\alpha|\le m$  and for some
$\beta$, with
$|\beta|=m-k+1$, the restriction of $D^\beta u$ to $U$ is not locally
 bounded.    
 
\section{Introduction}

  We will be concerned with local $C^\infty$ hypoellipticity in  
the following sense. A linear differential operator
operator $E$ on $\mathbb R^n$ is hypoelliptic if, whenever $u$ is a  
distribution such that the restriction of $Eu$ to an open
set $U\subset\mathbb R^n$ is in $C^\infty(U)$, then the restriction of  
$u$ to $U$ is also in $C^\infty(U)$. If
$E$ is hypoelliptic then it satisfies the following {\it a priori\/}  
estimates.
\begin{enumerate}
\item Given open sets $U, U'$ in $\mathbb R^n$ such that  
$U\subset\bar U\subset U'\subset\mathbb R^n$,  a nonnegative
integer $p$, and a real number $s_o$,   there exist  an integer $q$  
and a constant $C=C(U,p,q,s_o)$ such that
$$
\sum_{|\alpha|\le p} sup\,_{x\in U}|D^\alpha u(x)|\le  
C(\sum_{|\beta|\le q}sup\,_{x\in U'}|D^\beta Eu(x)|+\|u\|_{-s_o}),
$$
for all $u\in C_0^\infty(\mathbb R^n)$.
\item Given $\varrho,\varrho'\in C^\infty_0(\mathbb R^n)$ such that  
$\varrho'=1$ in a neighborhood of ${\rm supp}(\varrho)$, and
$s_o, s_1\in\mathbb R$,   there exist  $s_2\in\mathbb R$ and a  
constant $C=C(\varrho,\varrho',s_1,s_2,s_0)$ such that
$$
\|\varrho u\|_{s_1}\le C(\|\varrho'Eu\|_{s_2}+\|u\|_{-s_o}),
$$
for all $u\in C_0^\infty(\mathbb R^n)$.
\end{enumerate}
\noindent Assuming that $E$ is hypoelliptic and that $q$ is the  
smallest integer so that the first inequality above holds (for
large
$s_o$) then, if $q\le p$, we say that $E$  {\it gains} $p-q$  
derivatives in the sup norms and if $q\ge p$, we say that $E$
{\it loses} $q-p$ derivatives in the  sup norms. Similarly, assuming  
that $s_2$ is the smallest real number so that the
second inequality holds (for large $s_o$) then, if $s_2\le s_1$, we say  
that $E$ {\it gains} $s_1-s_2$ derivatives in the
Sobolev norms and if $s_2\ge s_1$, we say that $E$ {\it loses}  
$s_2-s_1$ derivatives in the Sobolev norms. In particular if
$E$ is of order $m$ and if $E$ is elliptic then $E$ gains exactly $m$  
derivatives in the Sobolev norms and gains exactly
$m-1$ derivatives in the sup norms. Here we will present hypoelliptic  
operators $E_k$ of order 2 which lose exactly $k-1$
derivatives in the Sobolev norms and lose at least $k$ derivatives in  
the sup norms.

  Loss of derivatives presents a very major difficulty: namely,  
how to derive the
{\it a priori\/} estimates? Such estimates depend on localizing the right-hand side and (because of the loss of derivatives)
the errors that arise are apparently always larger then the terms one  
wishes to estimate. This difficulty is overcome here
by the use of subelliptic multipliers in a microlocal setting. In this  
introduction I would like to indicate the ideas behind
these methods, which were originally devised to study hypoellipticity  
with gain of derivatives. It should be remarked that
that for global hypoellipticity the situation is entirely different; in  
that case loss of derivatives can occur and is well
understood but, of course, the localization problems do not arise.

 We will restrict ourselves to operators $E$ of second order  
of the form
$$
Eu=-\sum_{i,j}\frac{\partial}{\partial x_i}a_{ij}\frac{\partial  
u}{\partial x_j},
$$
where $(a_{ij})$ is a hermitian form with $C^\infty$ complex-valued  
components. If at some point $P\in\mathbb R^n$ the form
$(a_{ij}(P))$ has two nonzero eigenvalues of different signs then $E$  
is not hypoelliptic so that, without loss of generality,
we will assume that $(a_{ij})\ge0$.
\begin{defin} The operator $E$ is subelliptic at $P\in\mathbb R^n$ if  
there exists a neighborhood $U$ of $P$, a real
number $\varepsilon>0$, and a constant $C=(U,\varepsilon)$, such that
$$
\|u\|^2_\varepsilon\le C(|(Eu,u)|+\|u\|^2),
$$
for all $u\in C_0^\infty(U)$.
\end{defin}

Here the Sobolev norm $\|u\|_s$ is defined by
$$
\|u\|_s=\|\Lambda^su\|,
$$
and $\Lambda^su$ is defined by its Fourier transform, which is
$$
\widehat{\Lambda^su}(\xi)=(1+|\xi|^2)^{\frac{s}{2}}\hat u(\xi).
$$
We will denote by $H^s(\mathbb R^n)$ the completion of  
$C_0^\infty(\mathbb R^n)$ in the norm $\|\,\|_s$. If
$U\subset\mathbb R^n$ is open, we denote by $H^s_{\rm loc}(U)$ the set of  
all distributions on $U$ such that $\zeta
u\in H^s(\mathbb R^n)$ for all $\zeta\in C_0^\infty(U)$. The following  
result, which shows that subellipticity implies
hypoellipticity with a gain of $2\varepsilon$ derivatives in Sobolev  
norms, is proved in [KN].

\sdemo{Theorem} {\it Suppose that $E$ is subelliptic at each  
$P\in U\subset\mathbb R^n$. Then $E$ is hypoelliptic on
$U$. More precisely{\rm ,} if $u\in H^{-s_o}\cap H^s_{\rm loc}(U)$ and if $Eu\in  
H^s_{\rm loc}(U)${\rm ,} then} $u\in H^{s+2\varepsilon}_{\rm loc}(U)$.
\Enddemo

  In [K1] and [K2] I introduced subelliptic multipliers in  
order to  establish subelliptic estimates for the
$\bar\partial$-Neumann problem. In the case of
$E$, subelliptic multipliers are defined as follows.

\begin{defin} A {\it subelliptic multiplier} for $E$ at $P\in\mathbb  
R^n$ is a pseudodifferential operator $A$ of order zero,
defined on $C_0^\infty(U)$, where $U$ is a neighborhood of $P$, such  
that there exist  $\varepsilon>0$, and a constant
$C=C(\varepsilon,P,A)$, such that
$$
\|Au\|^2_\varepsilon\le C(|(Eu,u)|+\|u\|^2),
$$
for all $u\in C_0^\infty(U)$.
\end{defin}

If $A$ is a subelliptic multiplier and if $A'$ is a pseudodifferential  
operator whose principal symbol equals   the principal
symbol of $A$ then $A'$ is also a subelliptic multiplier. The existence  
of subelliptic estimates can be deduced from the
properties of the set symbols of subelliptic multipliers.  In the case of  
the $\bar\partial$-Neumann problem this leads to the
analysis of the condition of ``D'Angelo finite type.'' Catlin and  
D'Angelo, in [C] and [D'A], showed that D'Angelo
finite type is a necessary and sufficient condition for the  
subellipticity of the $\bar\partial$-Neumann problem. To
illustrate some of these ideas, in the case of an operator~$E$, we will  
recall H\"ormander's theorem on the sum of squares of
vector fields.

  Let $\{X_1,\dots,X_m\}$ be vector fields on a neighborhood of
the origin in $\mathbb R^n$.
\begin{defin} The vectorfields $\{X_1,\dots,X_m\}$ satisfy the {\it  
bracket condition} at the origin if the Lie algebra
generated by these vector fields evaluated at the origin is the tangent  
space.
\end{defin}

In [Ho], H\"ormander proved the following

\sdemo{Theorem} {\it If the vectorfields $\{X_1,\dots,X_m\}$ are  
real and if they satisfy the {\rm bracket condition}
at the origin then the operator $E=\sum X_j^2$ is hypoelliptic in a  
neighborhood of the origin.}
\Enddemo

The key point of the proof is to establish that for some  
neighborhoods of the origin $U$ there exist
$\varepsilon>0$ and $C=C(\varepsilon, U)$ such that
\begin{equation} \label{S}
\|u\|_\varepsilon^2\le C\big(\sum\|X_ju\|^2+\|u\|^2\big),
\end{equation}
for all $u\in C^\infty_0(U)$.
\medskip
Here is a brief outline of the proof of estimate (\ref{S}) using  
subelliptic multipliers. Note that
\begin{enumerate}
\item[1.] The operators $A_j=\Lambda^{-1}X_j$ are subelliptic multipliers  
with $\varepsilon=1$, that is
$$
\|A_ju\|_1^2\le C\big(\sum\|X_ju\|^2+\|u\|^2\big),
$$
for all $u\in C^\infty_0(U)$.
\item[2.] If $A$ is a subelliptic multiplier then $[X_j,A]$ is a  
subelliptic multiplier. (This is easily seen: we have
$X_j^*=-X_j+a_j$ since $X_j$ is real and
\begin{align*}
\|[X_j,A]u\|^2_{\frac{\varepsilon}{2}}&\le |(X_jAu,R^\varepsilon  
u)|+|(AX_ju,R^\varepsilon u)|\\
&\le |(Au,\tilde R^\varepsilon u)|+O(\|u\|^2)+|(Au,R^\varepsilon  
X_ju)|+|(AX_ju,R^\varepsilon u)|\\
&\le C\left(\|Au\|^2_\varepsilon+\sum\|X_ju\|^2+\|u\|^2\right),
\end{align*}
where $R^\varepsilon=\Lambda^\varepsilon[X_j,A]$ and $\tilde  
R^\varepsilon=[X^*_j,R^\varepsilon]$ are pseudodifferential
operators of order $\varepsilon$.)
\end{enumerate}
Now using the bracket condition and the above we see that 1 is a  
subelliptic multiplier and hence the estimate (\ref{S}) holds.
 
 The more general case, where the $a_{ij}$ are real but $E$  
cannot be expressed as a sum of squares
(modulo $L_2$) has been analyzed by Oleinik and Radkevic (see [OR]).  
Their result can also be obtained by use
of subelliptic multipliers and can then be connected to the geometric  
interpretation given by Fefferman and Phong in
[FP]. The next question, which has been studied fairly extensively, is  
what happens when subellipticity fails and yet there
is no loss. A striking example is the operator on $\mathbb R^2$ given by
$$
E=-\frac{\partial^2}{\partial x^2}-a^2(x)\frac{\partial^2}{\partial  
y^2},
$$
where $a(x)\ge0$ when $x\neq0$. This operator was studied by Fedii in  
[F], who showed that $E$ is always hypoelliptic, no matter
how fast $a(x)$ goes to zero as $x\to0$. Kusuoka and Stroock (see [KS])  
have shown that the operator on $\mathbb R^3$ given by
$$
E=-\frac{\partial^2}{\partial x^2}-a^2(x)\frac{\partial^2}{\partial  
y^2}-\frac{\partial^2}{\partial z^2},
$$
where $a(x)\ge0$ when $x\neq0$, is hypoelliptic if and only if  
$\lim_{x\to0}\log a(x)=0$.
Hypoellipticity when there is no loss but when the gain is smaller than  
in the subelliptic case has also been studied by Bell and
Mohamed [BM], Christ [Ch1], and  Morimoto [M]. Using subelliptic  
multipliers has provided new insights into these results (see
[K4]); for example Fedii's result is proved when $a^2$ is replaced by  
$a$ with the requirement that $a(x)>0$ when $x\neq0$. In
the case of the
$\bar\partial$-Neumann problem and of the operator $\Box_b$ on CR  
manifolds, subelliptic multipliers are used to established
hypoellipticity in certain situations where there is no loss of  
derivatives in Sobolev norms but in which the gain is weaker
than in the subelliptic case (see [K5]). Stein in [St] shows that the  
operator $\square_b+\mu$ on the Heisenberg group
$\mathcal H\subset\mathbb C^2$, with $\mu\neq0$, is analytic  
hypoelliptic but does not gain or lose any derivatives. In his thesis
Heller (see [He]),  using the methods developed by Stein in [St], shows  
that the fourth order operator $\square_b^2+X$ is
analytic hypoelliptic and that it loses derivatives (here $X$ denotes a  
``good'' direction).  In a recent work, C. Parenti and A. Parmeggiani studied classes of
pseudodifferential operators with large losses of derivatives (see [PP1]).

The study of subelliptic multipliers has led to the concept of  
multiplier ideal sheaves (see  [K2]). These have had
many applications notably Nadel's work on K\"ahler-Einstein metrics  
(see [N]) and numerous applications to algebraic
geometry. In algebraic geometry there are three areas in which  
multiplier ideals have made a decisive contribution: the
Fujita conjecture, the effective Matsusaka big theorem, and invariance  
of plurigenera; see, for example, Siu's article [S].
Up to now the use of subelliptic multipliers to study the  
$\bar\partial$-Neumann problem and the laplacian $\square_b$ has been
limited to dealing with Sobolev norms, Siu has developed a program to   
use multipliers for the $\bar\partial$-Neumann problem to
study H\"older estimates and to give an explicit construction of  the  
critical varieties that control the D'Angelo type. His
program leads to the study of the operator
$$
E=\sum_1^m X^*_jX_j,
$$
where the $\{X_1,\dots,X_m\}$ are {\it complex} vector fields  
satisfying the bracket condition. Thus Siu's program
gives rise to the question of whether the above operator $E$ is  
hypoelliptic and whether it satisfies the subelliptic estimate
(\ref{S}). These problems raised by Siu   have motivated my  
work on this paper. At first I found that if the bracket
condition involves only one bracket then (\ref{S}) holds with  
$\varepsilon=\frac{1}{4}$ (if the $X_j$ span without taking
brackets then $E$ is elliptic). Then I found a series of examples for  
which the bracket condition is satisfied with $k$
brackets, $k>1$, for which (\ref{S}) does not hold. Surprisingly I  
found that the operators in these examples
are hypoelliptic with a loss of $k-1$ derivatives in the Sobolev norms.  
  The method of proof involves calculations with
subelliptic multipliers and it seems very likely that it will be  
possible to treat the more general cases, that is when $E$
given by complex vectorfields and, more generally, when $(a_{ij})$ is  
nonnegative hermitian, along the same lines.

  The main results proved here are the following:

\sdemo{Theorem A} I{\it f $\{X_i,[X_i,X_j]\}$ span the  
complex
tangent space at the origin then a subelliptic estimate is satisfied{\it ,}
with
$\varepsilon =\frac{1}{2}$.}

\sdemo{Theorem B} {\it For $k\ge0$ there exist complex  
vector fields $X_{1k}$ and
$X_2$ on a neighborhood of the origin in $\mathbb R^3$ such that the  
two vectorfields $\{X_{1k}, X_2\}$ and their
commutators of order $k+1$ span the complexified tangent space at the  
origin{\rm ,} and when $k>0$ the subelliptic estimate
{\rm (\ref{S})} does not hold. Moreover{\rm ,} when $k>1${\rm ,} the operator  
$E_k=X_{1k}^*X_{1k}+X_2^*X_2$ loses $k$ derivatives in the sup
norms and $k-1$ derivatives in the Sobolev norms.}\Enddemo
 
Recently Christ (see [Ch2]) has shown that the operators  
$-\frac{\partial^2}{\partial s^2}+E_k$ on $\mathbb R^4$ are not
hypoelliptic when $k>0$.

\sdemo{Theorem C} {\it If $X_{1k}$ and $X_2$ are the  
vectorfields given in Theorem {\rm B} then the operator
$E_k=X_{1k}^*X_{1k}+X_2^*X_2$ is hypoelliptic. More precisely{\rm ,} if $u$  
is a distribution solution of $Eu=f$ with $u\in
H^{-s_0}(\mathbb R^3)$ and if
$U\subset\mathbb R^3$ is an open set such that $f\in H^{s_2}_{\rm loc}(U)${\rm ,}  
then $u\in H^{s_2-k+1}_{\rm loc}(U)$.}
\Enddemo
 
This paper originated with a problem posed by Yum-Tong Siu. The author  
wishes to thank Yum-Tong Siu and Michael Christ for
fruitful discussions of the material presented here.

\demo{Remarks}
In March 2005, after this paper had been accepted for publication, I circulated a preprint.
Then M. Derridj and D. Tartakoff proved analytic hypoellipticity for the operators constructed here (see [DT]). The work
of\break Derridj and Tartakoff used ``balanced'' cutoff functions to estimate the size of
derivatives starting with the
$C^\infty$ local hypoellipticity proved here; then Bove, Derridj, Tartakoff, and I (see [BDKT])
proved
$C^\infty$ local hypoellipticity using the balanced cutoff functions, starting from the estimates for functions with
compact support proved here. Also at this time, in [PP2],  Parenti and Parmeggiani, following their work in
[PP1],
 gave a different proof of hypoellipticity of the operators discussed here and in [Ch2].

\section{Proof of Theorem A}

The proof of Theorem A proceeds in the same way as given above in the  
outline of H\"ormander's theorem. It works only when one
bracket is involved because (unlike   the real case) $\bar X_j$ is not  
in the span of the $\{X_1,\dots,X_m\}$. The constant $\varepsilon=\frac{1}{2}$ is the largest possible,
since (as proved in [Ho]) this is already so when the $X_i$ are real.

   First note that $\|X^*_iu\|_{-\frac{1}{2}}^2\le
\|X_iu\|^2+C\|u\|^2$, since
\begin{eqnarray*}
\|X^*_iu\|_{-\frac{1}{2}}^2 &=&
(X^*_iu,\Lambda^{-1}X^*_iu)=(X^*_iu,P^0u)\\
&=&(u,X_iP^0 u)= -(u,P^0X_iu)+O(\|u\|^2);
\end{eqnarray*}
hence,
$$
\|X^*_iu\|_{-\frac{1}{2}}^2\le C\big(\sum\|X_ku\|^2+\|u\|^2\big),
$$
where $P^0=\Lambda^{-1}\bar{X}_i$ is a pseudodifferential operator of  
order zero. Then we have
\begin{eqnarray*}
\|X^*_iu\|^2 &=&
(u,X_iX^*_iu)=\|X_iu\|^2+(u,[X_i,X^*_i]u)\\
&=&\|X_iu\|^2+(\Lambda^{\frac{1}{2}}u,\Lambda^{-\frac{1}{2}}[X_i,X_i^*]u)\\
&\le& \|X_iu\|^2+C\|u\|^2_{\frac{1}{2}}.
\end{eqnarray*}

To estimate $\|u\|^2_{\frac{1}{2}}$ by
$C(\sum\|X_ku\|^2+\|u\|^2)$ we will estimate $\|Du\|^2_{-\frac{1}{2}}$ by\break
\vglue-10pt\noindent
$C(\sum\|X_ku\|^2+\|u\|^2)$ for all first order operators $D$. Thus it
suffices to estimate $Du$ when $D=X_i$ and when $D=[X_i,X_j]$. The estimate
is clearly satisfied if $D=X_i$, if $D=[X_i,X_j]$ we have
\begin{eqnarray*}
\|[X_i,X_j]u\|_{-\frac{1}{2}}^2
&=&(X_iX_ju,\Lambda^{-1}[X_i,X_j]u)-(X_jX_iu,\Lambda^{-1}[X_i,X_j]u)\\
&=&(X_iX_ju,P^0u)-(X_jX_iu,P^0u);
\end{eqnarray*}
the first term on the right is estimated by
\begin{eqnarray*}
(X_iX_ju,P^0u) &=&
(X_ju,X_i^*P^0u)=-(X_ju,P^0X^*_iu)+O(\|u\|^2+\|X_ju\|^2)\\ 
&\le & C(\|X_ju\|\|X^*_iu\|+\|u\|^2+\|X_ju\|^2)\\ &\le &
{\rm l.c.}\sum(\|X_ku\|^2+{\rm s.c.}\|X^*_iu\|^2+C\|u\|^2
\end{eqnarray*}
and the second term on the right is estimated similarly.
Combining these we have
$$
\|u\|^2_{\frac{1}{2}}\le  
C(\sum\|\frac{\partial u}{\partial x_i}\|^2_{-\frac{1}{2}}+\|u\|^2)
\le C(\sum\|X_ku\|^2+\|u\|^2)+{\rm s.c.}\|u\|^2_{\frac{1}{2}};
$$
hence
$$
\|u\|^2_{\frac{1}{2}}\le C\big(\sum\|X_ku\|^2+\|u\|^2\big)
$$
which concludes the proof of theorem A.

\section{The operators $E_k$}

  In this section we define the operators: $L,\bar L, X_{1k},  
X_2,$  and $E_k$.

Let $\mathfrak H$ be the hypersurface in $\mathbb C^2$
given by:
$$
\Re(z_2)=-|z_1|^2.
$$
We identify $\mathbb R^3$ with the Heisenberg group represented by  
$\mathfrak H$ using the mapping
$\mathfrak H\to\mathbb R^3$ given by $x=\Re
z_1,\, y=\Im z_1,\, t=\Im z_2$. Let $z=x+\sqrt{-1}\,y$. Let
$$
L=\frac{\partial}{\partial z_1}-2\bar z_1\frac{\partial}{\partial z_2}=
\frac{\partial}{\partial z}+\sqrt{-1}\bar z\frac{\partial}{\partial t}
$$
and
$$
\bar L=\frac{\partial}{\partial\bar  
z_1}-2z_1\frac{\partial}{\partial\bar z_2}=
\frac{\partial}{\partial\bar z}-\sqrt{-1}\,z\frac{\partial}{\partial t}.
$$
Let $X_{1k}$ and $X_2$ be the restrictions to $\mathfrak H$ of the  
operators
  $$
X_{1k}=\bar z_1^kL=\bar z^k\frac{\partial}{\partial z}+\sqrt{-1}\,\bar  
z^{k+1}\frac{\partial}{\partial t}.
$$
We set
$$
X_2=\bar L= \frac{\partial}{\partial\bar  
z}-\sqrt{-1}\,z\frac{\partial}{\partial t}
$$
and
$$
E_k=X^*_{1k}X_{1k}+X^*_2X_2=-\bar L|z|^{2k}L-L\bar L.
$$

By induction on $j$ we define the commutators $A^j_k$ setting  
$A^1_k=[X_{1k},X_2]$ and
$A^j_k=[A^{j-1}_k,X_2]$. Note that  
$X_2, A^k_k$ and $A^{k+1}_k$ span the tangent space of~$\mathbb R^3$.

\section{Loss of derivatives (part I)}

In this section we prove that the subelliptic
estimate (\ref{S}) does not hold when $k\ge1$. We also prove a  
proposition  which gives the loss of derivatives in the sup
norms which is part of Theorem B. To complete the proof of Theorem B,  
by  establishing loss in the Sobolev norms, we will
use additional microlocal analysis of $E_k$, the proof of Theorem B is  
completed in Section 6.

\begin{defin} If $U$ is a neighborhood of the origin then $\varrho\in  
C_0^\infty(U)$ is real-valued and is defined as follows
$\varrho(z,t)=\eta(z)\tau(t)$, where $\eta\in C_0^\infty(\{z\in\mathbb  
C\mid|z|<2\})$ with $\eta(z)=1$ when $|z|\le 1$ and
$\tau\in C_0^\infty(\{t\in\mathbb R\mid|t|<2a\})$ with $\tau(t)=1$ when  
$|t|\le a$.
\end{defin}

 The following proposition shows that the subelliptic estimate  
(\ref{S}) does not hold when $k>0$.

\begin{prop} If $k\ge1$ and if there exist  a neighborhood $U$ of the  
origin and constants $s$ and $C$ such that
$$
\|u\|^2_s\le C(\|\bar z^kLu\|^2+\|\bar Lu\|^2),
$$
for all $u\in C_0^\infty(U)${\rm ,} then $s\le0$.
\end{prop}

  {\it Proof.} Let $\lambda_0$ and $a$ be sufficiently large so  
that the support of $\varrho(\lambda z,t)$ lies in
$U$ when $\lambda\ge\lambda_0$. We define $g_\lambda$ by
$$
g_\lambda(z,t)=\varrho(\lambda  
z,t)\exp(-\lambda^{\frac{5}{2}}(|z|^2-it)).
$$
Note that $L\eta(z)=\bar L\eta(z)=0$ when $|z|\le1$, that  
$L(\tau)=i\bar z\tau'$, and that $\bar L(\tau)=-iz\tau'$.
Setting $R^\lambda v(z,t)=v(\lambda z,t)$, we have:
\begin{align} \label{Lg}
\bar z^kL(g_\lambda)&=(\lambda\bar z^k(R^\lambda L\eta)\tau+i\bar  
z(R^\lambda\eta)\tau' +\lambda^\frac{5}{2}\bar
zR^\lambda\varrho)\exp(-\lambda^\frac{5}{2}(|z|^2+it))
\end{align}
and
\begin{equation} \label{barLg}
\bar L(g_\lambda)=(\lambda (R^\lambda\bar  
L\eta)\tau-iz(R^\lambda\eta)\tau')\exp(-\lambda^\frac{5}{2}(|z|^2+it)).
\end{equation}
Note that the restriction of $|g_\lambda|$ to $\mathfrak H$ is
$$
|g_\lambda(z,t)|=\varrho(\lambda z,t)\exp(-\lambda^\frac{5}{2}|z|^2).
$$
Now we have, using the changes of variables: first $(z,t)\mapsto  
(\lambda^{-1}z,t)$ and then
$z\mapsto \lambda^{-\frac{1}{4}}z$
\begin{eqnarray*}
\|g_\lambda\|^2&=&\frac{C}{\lambda^2}\int_{\mathbb  
R^2}\eta(z)^2\exp(-2\lambda^\frac{1}{2}|z|^2) dxdy\\
&\ge&\frac{C}{\lambda^2}\int_{\mathbb  
R^2}\exp(-2\lambda^\frac{1}{2}|z|^2)dxdy
-\frac{C}{\lambda^2}\int_{|z|\ge1}\exp(-2\lambda^\frac{1}{2}|z|^2)dxdy\\
&\ge&\frac{C}{\lambda^\frac{5}{2}}-
\frac{C}{\lambda^2}\exp(-\lambda^{\frac{1}{2}})\int_{\mathbb  
R^2}\exp(-\lambda^{\frac{1}{2}}|z|^2)dxdy\\
&\ge&\frac{C}{\lambda^{\frac{5}{2}}}- 
\frac{C}{\lambda^{\frac{5}{2}}}\exp(-\lambda^{\frac{1}{2}}).
\end{eqnarray*}
 Then we have
$$
\|g_\lambda\|^2\ge  
\frac{{\rm const}.}{\lambda^{\frac{5}{2}}}
$$
for sufficiently large $\lambda$. Further, using the above coordinate  
changes to  estimate the individual terms in (\ref{Lg}) and
in (\ref{barLg}), we have
\begin{align*}
&\|\bar z^k\lambda (R^\lambda  
L\eta)\tau\exp(-\lambda^{\frac{5}{2}}(|z|^2-it)\|^2+
\|\lambda (R^\lambda\bar  
L\eta)\tau\exp(-\lambda^{\frac{5}{2}}(|z|^2-it)\|^2\\
&\qquad\le  
C\exp(-\lambda^{\frac{1}{2}})\int_{|z|\ge1}\exp(- 
\lambda^{\frac{1}{2}}|z|^2)dxdy\le
\frac{C}{\lambda^{\frac{1}{2}}}\exp(-\lambda^{\frac{1}{2}}),\end{align*}
\begin{eqnarray*}
&&\hskip-36pt \||z|(R^\lambda\eta)\tau')\exp(-\lambda^\frac{5}{2}(|z|^2+it))\|^2\\
&&\qquad\le
\frac{C}{\lambda^2}\int|z|^2\exp(-2\lambda^{\frac{1}{2}}|z|^2)dxdy 
 \le \frac{C}{\lambda^4},\\
\noalign{\noindent \rm and}
&&\hskip-36pt
\|\lambda^\frac{5}{2}\bar  
z^{k+1}R^\lambda\varrho)\exp(-\lambda^\frac{5}{2}(|z|^2+it))\|^2\\
&&\qquad\le
C\lambda^{1-2k}\int|z|^{2k+2}\exp(-2\lambda^{\frac{1}{2}}|z|^2)dxdy\le
\frac{C}{\lambda^{\frac{5k}{2}}}.
\end{eqnarray*}
Hence, if $k\ge1$, we have
$$
\|\bar z^kLg_\lambda\|^2+\|\bar Lg_\lambda\|^2\le  
\frac{C}{\lambda^{\frac{5}{2}}}.
$$
Since $|x|\le\frac{2}{\lambda}$ on the support of $g_\lambda$ then we  
conclude, from the lemma proved below, that given
$\varepsilon$ there there exists $C$ such that
$$
\lambda^\varepsilon\|g_\lambda\|\le C\|g_\lambda\|_\varepsilon,
$$
for sufficiently large $\lambda$. It the follows that, if $k\ge1$ then  
the subelliptic estimate
$$
\|g_\lambda\|_\varepsilon^2\le C(\|\bar z^kLg_\lambda\|^2+\|\bar  
Lg_\lambda\|^2)
$$
implies that $\lambda^{2\varepsilon-\frac{5}{2}}\le  
C\lambda^{-\frac{5}{2}}$ which is a contradiction and thus the  
proposition
follows. The following lemma then completes the proof. For completeness
we include a proof which is along the lines given in [ChK].
\bigskip
\begin{lemma}
Let $Q_\delta$ denote a bounded open set contained in the \/{\rm ``}\/slab\/{\rm ''}\break
$\{x\in\mathbb R^n\mid|x_1|\le\delta\}$. Then{\rm ,}
for each
$\varepsilon>0${\rm ,} there exists
$C=C(\varepsilon)>0$ such that
\begin{equation} \label{CS}
\|u\|\le C\delta^\varepsilon\|u\|_\varepsilon,
\end{equation}
for all $u\in C^\infty_0(Q_\delta)$ and $\delta>0$.
\end{lemma}

\Proof  Note that the general case follows from the case
of $n=1$. Since, writing $x=(x_1,x')$, if for each fixed
$x'$ we have $\|u(\cdot,x')\|\le
C\delta^\varepsilon\|u(\cdot,x')\|_\varepsilon$ then, after integrating
with respect to
$x'$ and noting that\\
$(1+\xi_1^2)^\varepsilon\le(1+|\xi|^2)^\varepsilon$, we obtain the
desired estimate. So we will assume
that $n=1$ and  set $x=x_1$  and $\xi=\xi_1$. We define $\||u\||_s$ by
$$
\||u\||_s^2=\int|\xi|^{2s}|\hat u(\xi)|^2d\xi.
$$
We will show that, if $s\ge0$, there exists a constant $C$ such that
$$
\|u\|_s\le C\||u\||_s,
$$
for all $u\in C_0^\infty((-1,1))$.
First we have
$$
|\hat u(\xi)|=|\int e^{-ix\xi}u(x)dx|\le \sqrt{2}\|u\|.
$$
Next, if $|\xi|\le a\le1$,  
$$
(1+\xi^2)^s|\hat u(\xi)|^2\le2^{s+1}\|u\|^2
$$
and
\begin{multline*}
\int_{-\infty}^{\infty}(1+\xi^2)^s|\hat u(\xi)|^2d\xi\\=\int_{|\xi|\le  
a}\dots
+\int_{|\xi|>a}\dots
\le2^{s+2}a\|u\|^2+\left(\frac{1}{a^2}+1\right)^s\||u\||_s^2.
\end{multline*}
Hence if $a$ is small we obtain $\|u\|_s\le C\||u\||_s$, as required.   
If ${\rm supp}(u)\subset(-\delta,\delta)$ then
set $u_\delta(x)=u(\delta x)$ so that ${\rm supp}( u_\delta)\subset(-1,1)$.  
Now
$$
\hat{u_\delta}(\xi)=\frac{1}{\delta}\hat u\left(\frac{\xi}{\delta}\right)
$$
so that $\|u_\delta\|^2=\frac{1}{\delta}\|u\|^2$ and  
$\||u_\delta\||^2_s=\delta^{2s-1}\||u\||^2_s$ which concludes the proof.
\Enddemo

 Next we prove that $E_k$ loses at least $k$ derivatives in  
the sup norms.

\begin{prop} If for some open sets $U$ and $U'${\rm ,} with $\bar U\subset  
U'${\rm ,} and for each $s_0$ there exists a constant $C=C(s_0)$
such that
\begin{equation} \label{sup}
\sum_{|\alpha|\le p}\sup_{x\in U}|D^\alpha u(x)|\le C\Big(\sum_{|\beta|\le  
q}\sup_{x\in U'}|D^\beta E_ku(x)|
+\|u\|_{-s_0}\Big),
\end{equation}
for all $u\in C_0^\infty(\mathbb R^3)${\rm ,} then $q\ge p+k$.
\end{prop}

\Proof  If $\delta>0$ define $u_\delta$ by
$$
u_\delta=(|z|^2-\sqrt{-1}t)^p\log(|z|^2+\delta-\sqrt{-1}t),
$$
where $\log$ denotes the branch of the logarithm that takes reals into  
reals. Since $u_\delta$ is the restriction of
$(-z_2)^p\log(-z_2+\delta)$ to $\mathfrak H$ we have $\bar  
Lu_\delta=0$. Then we  have
$$
{\rm lim}_{\delta\to0}|D^p_tu_\delta(0)|=\infty.
$$
Further
\begin{align*}
E_ku_\delta&=-\bar  
L|z|^{2k}Lu_\delta\\
&=2k|z_1|^{2k}\!\!\left(-p(-z_2)^{p-1}\log(-z_2+\delta)+(- 
z_2)^p\log(-z_2+\delta)+
\frac{(-z_2)^p}{(-z_2+\delta)}\right)\\
&=2k|z|^{2k}\!\!\left(p(|z|^2-\sqrt{-1}t)^{p-1}\log(|z|^2+\delta-\sqrt{- 
1}t)+
\frac{(|z|^2-\sqrt{-1}t)^p}{|z|^2+\delta-\sqrt{-1}t}\right).
\end{align*} \pagebreak

\noindent
Note that $\|u_\delta\|_{-s_0}$ is bounded independently of $\delta$  
when $s_0\ge3$.
Thus, when $q\le p+k-1$, we have
$$
\sum_{|\beta|\le q}\sup_{x\in U'}|D^\beta E_ku_\delta(x)|\le {\rm const}.
$$
with the constant independent of $\delta$. Hence, applying (\ref{sup})  
to $u_\delta$ we obtain
$q\ge p+k$. This concludes the proof of the proposition.

\section{Notation}

In this section we set down some notation which will be used throughout  
the rest of the paper.
\begin{enumerate}
\item[1.] Associated to the cutoff function $\varrho$ defined in Definition  
1, is a $C^\infty$ function $\mu$ such that
$L\varrho=\bar z\mu$ and $\bar L\varrho=z\bar\mu$ (Such a $\mu$ exists  
since
$$
L\varrho(z,t)=D_z\eta(z)\tau(t)+i\bar z\eta(z)D_t\tau(t).
$$
Since $D_z\eta(z)=0$ in a neighborhood of $z=0$ we can set  
$\mu(z,t)=\frac{D_z\eta(z)}{\bar z}\tau(t)+i\eta(z)D_t\tau(t)$.)
\item[2.] Given cutoff functions $\varrho, \varrho'$, as in Definition 1,  
with $\varrho'=1$ in a neighborhood of the support of
$\varrho$, then we denote by $\{\varrho_i\}$ a special sequence of  
cutoff functions, each of which satisfies Definition 1 and
such that: $\varrho_1=\varrho$, $\varrho'=1$ in a neighborhood of  
$\bigcup\varrho_i$, and $\varrho_{i+1}=1$ in a neighborhood of
the support of  $\varrho_i$.
\item[3.] The abbreviations ``${\rm s.c.}$'' and ``${\rm l.c.}$'' will be used for  
``small constant'' and ``large constant'', respectively in the
following sense. $\mathcal Au\le {\rm s.c.}\mathcal Bu+{\rm l.c.}\mathcal Cu$ means  
that given any constant ${\rm s.c.}$ there exists a constant
${\rm l.c.}$ such that the inequality holds for all $u$ in some specified  
class.
\item[4.] We will use $\|u\|_{-\infty}$ to denote the following. Given  
$\mathcal Au$, the expression $\mathcal Au\le
\|u\|_{-\infty}$ means that: if for any $s_o$ there exist a constant  
$C=C(s_o)$ such that $\mathcal Au\le C\|u\|_{-s_o}$
holds for all $u$ in some specified class.
\end{enumerate}

\section{Microlocalization on the Heisenberg group}

Denote by $T$ the vector field defined by
$$
T=\frac{1}{\sqrt{-1}}\frac{\partial}{\partial t}.
$$
Then 
$$
[L,\bar L]=[\frac{\partial}{\partial z}+\sqrt{-1}\,\bar  
z\frac{\partial}{\partial t},\,\frac{\partial}{\partial\bar
z}-\sqrt{-1}\,z\frac{\partial}{\partial t}]=2T.
$$
The following simple formula, which is obtained by integration by  
parts, is the starting point of all the estimates connected
with the operators $E_k$.

\begin{lemma} For $u\in C^\infty_0(\mathbb R^3)$ we have
\begin{equation} \label{basic}
\|Lu\|^2=2(Tu,u)+\|\bar Lu\|^2.
\end{equation}
\end{lemma}
\vglue8pt

\Proof  Since $L^*=-\bar L$ and $\bar L^*=-L$, we have
$$
\|Lu\|^2=(Lu,Lu)=-(\bar LLu,u)=-([\bar L,L]u,u)-(L\bar  
Lu,u)=2(Tu,u)+\|\bar Lu\|^2.
$$
We set $x_1=\Re z, x_2=\Im z,$ and $x_3=t$ and denote the dual  
coordinates by $\xi_1,\xi_2,$ and $\xi_3$.
For $(\alpha,t_0)\in\mathbb C\times\mathbb R$ we define
$$
z^\alpha=z-\alpha \ \ {\rm and}\ \  
x_3^\alpha=-2\alpha_2x_1+2\alpha_1x_2+x_3-t_0,
$$
where $\alpha_1=\Re\alpha$ and $\alpha_2=\Im\alpha$. Then
$$
L=\frac{\partial}{\partial z^\alpha}+i\bar  
z^\alpha\frac{\partial}{\partial x_3^\alpha}
$$
and
$$
\bar L=\frac{\partial}{\partial\bar  
z^\alpha}-iz^\alpha\frac{\partial}{\partial x_3^\alpha}.
$$
We set $x_1^\alpha=x_1-\alpha_1$, $x_2^\alpha=x_2-\alpha_2$, and  
$x^\alpha=(x^\alpha_1,x^\alpha_2,x^\alpha_3)$. Let $\mathcal
F_\alpha$ denote the the Fourier transform in the
$x^\alpha_j$ coordinates; that is
$$
\mathcal F_\alpha u(\xi)=\int  
e^{-ix^\alpha\cdot\xi}u(x^\alpha)dx^\alpha_1x^\alpha_2x^\alpha_3.
$$

\begin{defin} Let $S^2=\{\xi\in\mathbb R^3\mid |\xi|=1\}$ be the unit  
sphere.
Suppose that $\mathcal U, \mathcal U_1$ are open subsets of $S^2$ with  
$\bar{\mathcal U}_1\subset\mathcal U$.
For each such pair of open sets we define a set of $\gamma\in  
C^\infty(\mathbb R^3)$, with $\gamma\ge0$, such that
\begin{enumerate}
\item[1.] $\gamma(\frac{\xi}{|\xi|})=\gamma(\xi)$ when $|\xi|\ge1$.
\item[2.] $\gamma(\xi)=1$ when $\xi\in\mathcal U_1$.
\item[3.] $\gamma(\xi)=0$ when $\xi\in S^2-\mathcal U$.
\end{enumerate}
To such a $\gamma$ and $\alpha\in\mathbb C$ we associate the operator  
$\Gamma_\alpha$ defined by
$$
\mathcal F_\alpha\Gamma_\alpha u(\xi)=\gamma(\xi)\mathcal F_\alpha  
u(\xi).
$$
Let $\mathcal U^+,\mathcal U_1^+,
\mathcal U^0,\mathcal U^0_1,\mathcal U^-$, and $\mathcal U^-_1$ be open  
subsets of $S^2$ defined as follows.
\begin{align*}
&\mathcal U^+=\Big\{\xi\in S^2\mid \xi_3>\frac{5}{9}\Big\},\ \ \ \mathcal  
U^+_1=\Big\{\xi\in S^2\mid \xi_3>\frac{4}{9}\Big\},\\
&\mathcal U^0=\Big\{\xi\in S^2\mid |\xi_3|<\frac{5}{6}\Big\},\ \ \ \mathcal  
U^0_1=\Big\{\xi\in S^2\mid |\xi_3|<\frac{2}{3}\Big\},\\*
&\mathcal U^-=\{\xi\in S^2\mid -\xi\in\mathcal U^+\},\ {\rm and}\ \  
\mathcal U^-_1= \{\xi\in S^2\mid -\xi\in\mathcal U^+_1 \}.
\end{align*}
We denote by $\gamma^+$, $\gamma^0$, and $\gamma^-$ the corresponding  
functions and require further that
$\gamma^+(\xi)=\gamma^-(\xi)=0$ when $|\xi|\le\frac{1}{2}$ and  
$\gamma^0(\xi)=1$ when $\frac{\xi}{|\xi|}\in\mathcal U^0_1$. The
sets of these functions will be denoted by $\mathcal G^+$, $\mathcal  
G^0$, and $\mathcal G^-$, respectively.  The corresponding
operators are denoted by $\Gamma^+_\alpha$, $\Gamma^0_\alpha$, and  
$\Gamma^-_\alpha$. The sets of these operators will be
denoted by $\mathfrak G^+_\alpha$, $\mathfrak G^0_\alpha$, and  
$\mathfrak G^-_\alpha$, respectively. Given
$(\alpha,t_0)\in\mathbb C\times\mathbb R$ the functions  
$\Gamma^+_\alpha u$, $\Gamma^0_\alpha u$, and $\Gamma^-_\alpha u$ will
be referred to as microlocalizations of $u$ at $(\alpha,t_0)$ in the  
regions $+$, $0$, and $-$, respectively.
\end{defin}

The following lemma shows that the $0$ microlocalization is elliptic  
for the operators $L$ and $\bar L$.
In our estimates we will often encounter error terms which can be  
bounded by $C_{s_0}\|u\|_{-s_0}$ for every
$s_0$; abusing notation we will bound such terms by  
$``\|u\|_{-\infty}"$.

\begin{lemma} If $U$ is a neighborhood of $(\alpha,t_0)$ and if  
$\gamma^0,\tilde\gamma^0\in\mathcal G^0$ with
$\tilde\gamma^0=1$ in a neighborhood of the support of $\gamma^0$ then  
there exist constants $a>0$ and
$C>0$ such that{\rm ,} if $|z-\alpha|<a$  on $U${\rm ,} then
$$
\|\Gamma^0_\alpha u\|_1\le C(\|\Gamma^0_\alpha  
Lu\|+\|\tilde\Gamma^0_\alpha u\|+\|u\|_{-\infty})
$$
and
$$
\|\Gamma^0_\alpha u\|_1\le C(\|\Gamma^0_\alpha\bar  
Lu\|+\|\tilde\Gamma^0_\alpha u\|+\|u\|_{-\infty}),
$$
for all $u\in C^\infty_0(U)$.
\end{lemma}

\Proof  If $\xi\in\mathcal U^0$ and if $|\xi|\ge1$ then  
$|\xi_3|\le\frac{5}{6}|\xi|$. Thus, if $\xi\in\mathcal U^0$,
then $|\xi|\le6(|\xi_1|+|\xi_2|)+1$. Now,
$$
\|\Gamma^0_\alpha u\|^2_1\le C\Big(\sum_1^2\|\frac{\partial}{\partial  
x^\alpha_j}\Gamma^0_\alpha u\|^2+
\|\tilde\Gamma^0_\alpha u\|^2+\|u\|_{-\infty}^2\Big).
$$
Let $U'\supset\bar U$ be an open set such that $|z-\alpha|>2a$ on $U'$  
and let $\varphi\in C_0^\infty(U')$ satisfying
$\varphi=1$ in a neighborhood of $\bar U$. Then
\begin{align*}
\|\Gamma^0_\alpha u\|^2_1&\le C\Big(\sum_1^2\|\frac{\partial}{\partial  
x^\alpha_j}\Gamma^0_\alpha\varphi u\|^2+
\|\tilde\Gamma^0_\alpha u\|^2+\|u\|_{-\infty}^2\Big)\\
&\le C'\Big(\sum_1^2\|\frac{\partial}{\partial  
x^\alpha_j}\varphi\Gamma^0_\alpha u\|^2+
\|\tilde\Gamma^0_\alpha u\|^2+\|u\|_{-\infty}^2\Big)\\
&\le C''(\|L\varphi\Gamma^0_\alpha u\|^2+\|\bar L\varphi\Gamma^0_\alpha  
u\|^2\\
& \phantom{\le C''(} +
\max_{U'}|z-\alpha|^2\|\frac{\partial}{\partial  
x^\alpha_3}\Gamma^0_\alpha u\|^2+
\|\tilde\Gamma^0_\alpha u\|^2+\|u\|_{-\infty}^2)\\
&\le C''(\|\Gamma^0_\alpha Lu\|^2+\|\Gamma^0_\alpha\bar Lu\|^2+
4a^2\|\Gamma^0_\alpha u\|_1^2+\|\tilde\Gamma^0_\alpha  
u\|^2+\|u\|_{-\infty}^2).
\end{align*}
Hence, taking $a$ suitably small we obtain
$$
\|\Gamma^0_\alpha u\|^2_1\le C(\|\Gamma^0_\alpha  
Lu\|^2+\|\Gamma^0_\alpha\bar L u\|^2+
\|\tilde\Gamma^0_\alpha u\|^2+\|u\|_{-\infty}^2).
$$
Furthermore, substituting $\varphi\Gamma^0_\alpha u$ for $u$ in  
(\ref{basic}), we have
\begin{align*}
\|L\varphi\Gamma^0_\alpha u\|^2&=2(T\varphi\Gamma^0_\alpha  
u,\varphi\Gamma^0_\alpha u)+\|\bar L\varphi\Gamma^0_\alpha u\|^2\\
&\le {\rm s.c.}\|\frac{\partial}{\partial x^\alpha_3}\Gamma^0_\alpha u\|^2+
{\rm l.c.}(\|\tilde\Gamma^0_\alpha  
u\|^2+\|u\|_{-\infty}^2)+C\|\Gamma^0_\alpha\bar Lu\|^2\\
&\le {\rm s.c.}\|\Gamma^0_\alpha u\|^2_1+{\rm l.c.}(\|\tilde\Gamma^0_\alpha  
u\|^2+\|u\|_{-\infty}^2)+C\|\Gamma^0_\alpha\bar Lu\|^2,
\end{align*}
and since
$$
\|L\varphi\Gamma^0_\alpha u\|^2\le C(\|\Gamma^0_\alpha  
Lu\|^2+\|\tilde\Gamma^0_\alpha u\|^2+\|u\|_{-\infty}^2)
$$
we get
$$
\|\Gamma^0_\alpha u\|_1\le C(\|\Gamma^0_\alpha\bar  
Lu\|+\|\tilde\Gamma^0_\alpha u\|^2+\|u\|_{-\infty}^2).
$$
Similarly  we obtain
$$
\|\Gamma^0_\alpha u\|_1\le C(\|\Gamma^0_\alpha  
Lu\|+\|\tilde\Gamma^0_\alpha u\|^2+\|u\|_{-\infty}^2).
$$
This completes the proof of the lemma.

\begin{lemma}If $R^s$ is a pseudodifferential operator of order $s$  
then there exists $C$ such that
\begin{eqnarray*}
\|[R^s,\Gamma^+_\alpha]u\|&\le &C(\|\Gamma^0_\alpha  
u\|_{s-1}+\|u\|_{-\infty})\\
\noalign{\noindent and} 
|[R^s,\Gamma^-_\alpha]u\|&\le &C(\|\Gamma^0_\alpha  
u\|_{s-1}+\|u\|_{-\infty}).
\end{eqnarray*}
\end{lemma}
\vskip8pt

\Proof  Since $\gamma^0=1$ on a neighborhood of the support of the  
derivatives of $\gamma^+$ it also equals one on a
neighborhood of the support of  the symbol of $[R^s,\Gamma^+_\alpha]$.  
Hence
$[R^s,\Gamma^+_\alpha]=[R^s,\Gamma^+_\alpha]\Gamma^0_\alpha+R^{- 
\infty}$, where
$R^{-\infty}$ is a pseudodifferential operator whose symbol is  
identically zero. The same argument works for the term
$[R^s,\Gamma^-_\alpha]$ and the lemma follows.
 
\begin{defin} For each $s\in\mathbb R$ we define the operator  
$\Psi^s_\alpha$ as follows. Let $\mathcal U^*$ and
$\mathcal U^*_1$ be open sets in $S^2$ such that $\mathcal U^*=\{\xi\in  
S^2\mid |\xi_3|>\frac{1}{6}$ and
$\mathcal U_1^*=\break\{\xi\in S^2\mid |\xi_3|>\frac{1}{3}$. Let $\gamma^*$  
be the function on $R^3$ associated with
$\mathcal U^*,\ \mathcal U^*$ such that $\gamma^*(\xi)=0$ when  
$|\xi|\le\frac{1}{3}$ and $\gamma^*(\xi)=1$ in the region
$\{\xi\in\mathbb R^3 \mid \frac{\xi}{|\xi|}\in\mathcal U^*_1\ {\rm  
and}\ |\xi|\ge\frac{1}{2}\}$. Then we set
$\psi^s(\xi)=(1+|\xi_3|^2)^{\frac{s}{2}}\gamma^*(\xi)$ and define  
$\Psi^s_\alpha$ by
$$
\mathcal F_\alpha\Psi^s_\alpha u(\xi)=\psi^s(\xi)\mathcal F_\alpha  
u(\xi).
$$
\end{defin}

Note that there exist positive constants $c$ and $C$ such that
$$
c(1+|\xi|^2)^{\frac{s}{2}}\gamma^*(\xi)\le\psi^s(\xi)\le  
C(1+|\xi|^2)^{\frac{s}{2}}\gamma^*(\xi).
$$
Hence $\|\Psi^s_\alpha\Gamma^+_\alpha u\|\sim\|\Gamma^+_\alpha u\|_s$
and $\|\Psi^s_\alpha\Gamma^-_\alpha u\|\sim\|\Gamma^-_\alpha u\|_s$; by  
$\sim$ we mean that they differ by an operator of
order $-\infty$. Also, since $\gamma^*=1$ on the supports of $\gamma^+$  
and $\gamma^-$, we have
$$
\Psi^s\Psi^{s'}\Gamma^+_\alpha\sim\Psi^{s+s'}\Gamma^+_\alpha\ {\rm  
and}\ \Psi^s\Psi^{s'}\Gamma^-_\alpha\sim
\Psi^{s+s'}\Gamma^-_\alpha.
$$

\begin{lemma}There exists $C$ such that
$$
\|\Gamma^+_\alpha\bar Lu\|^2+\|\Gamma^+_\alpha u\|_{\frac{1}{2}}^2\le
C(\|\Gamma^+_\alpha Lu\|^2+\tilde\Gamma^+_\alpha  
u\|^2+\|u\|_{-\infty}^2),
$$
and
$$
\|\Gamma^-_\alpha Lu\|^2+\|\Gamma^-_\alpha u\|_{\frac{1}{2}}^2\le
C(\|\Gamma^-_\alpha\bar Lu\|^2+\tilde\Gamma^-_\alpha  
u\|^2+\|u\|_{-\infty}^2),
$$
for all $u\in C_0^\infty(U)$.
\end{lemma}
\Proof  Taking $\varphi\in C_0^\infty$ with $\varphi=1$ in a  
neighborhood of $\bar U$ we substitute
$\varphi\Gamma^+_\alpha u$ for $u$ in (\ref{basic}) and obtain
$$
\|L\varphi\Gamma^+_\alpha u\|^2=2(T\varphi\Gamma^+_\alpha  
u,\varphi\Gamma^+_\alpha u)+
\|\bar L\varphi\Gamma^+_\alpha u\|^2.
$$
Now, we have
$$
(T\varphi\Gamma^+_\alpha u,\varphi\Gamma^+_\alpha  
u)=(T\Gamma^+_\alpha\varphi u,\varphi\Gamma^+_\alpha u)+
O(\tilde\Gamma^+_\alpha u\|^2+\|u\|_{-\infty}^2).
$$
Since $\mathcal F_\alpha(Tu)=\xi_3\mathcal F_\alpha(u)$ we have  
$T\Gamma^+_\alpha\sim\Psi^1\Gamma^+_\alpha\sim
\Psi^{\frac{1}{2}}\Psi^{\frac{1}{2}}\Gamma^+_\alpha$ and
$$
(T\varphi\Gamma^+_\alpha u,\varphi\Gamma^+_\alpha u)=
\|\Psi^{\frac{1}{2}}_\alpha\Gamma^+_\alpha  
u\|^2+O(\tilde\Gamma^+_\alpha u\|^2+\|u\|_{-\infty}^2).
$$
This proves the first part of the lemma, the second follows from the  
fact that $|\xi_3|\gamma^-(\xi)=-\xi_3\gamma^+(-\xi)$.
Then  
$\Psi^1\Gamma^-_\alpha\sim\Psi^{\frac{1}{2}}\Psi^{\frac{1}{2}}\Gamma^- 
_\alpha$, thus concluding the proof.

\section{Loss of derivatives (part II) \\ Conclusion of the proof of  
Theorem B}

  In this section we conclude the proof of Theorem B by showing  
that if $k\ge2$ then $E_k$ loses at least $k-1$
derivatives in the Sobolev norms.

\begin{prop} Suppose that there exist two neighborhoods of the origin  
$U$ and $U'${\rm ,} with $\bar U\subset U'${\rm ,} and
real  numbers $s_1$ and $s_2$ such that if $\varrho,\varrho'\in  
C_0^\infty(U')$ with $\varrho=1$ on $U$ and
$\varrho'=1$ in a neighborhood of the support of $\varrho${\rm ,} and if for  
any real number $s_0$ there exists a constant
$C=C(\varrho,\varrho',s_0)$ such that
\begin{equation} \label{hypest}
\|\varrho u\|_{s_1}\le C(\|\varrho' E_ku\|_{s_2}+\|u\|_{-s_0}),
\end{equation}
for all $u\in \mathcal S$, then $s_2\ge s_1+k-1$. Here $\mathcal S$  
denotes the Schwartz space of rapidly decreasing functions.
\end{prop}

\Proof  Let $\{\varrho_i\}$ and $\{\varrho'_i\}$ be sequences of  
cutoff functions in $C^\infty_0(U)$ and $C^\infty_0(U')$,
respectively. We assume that $\varrho_(z,t)=\eta_i(|z|)\tau_i(t)$ and  
$\varrho'_i(z,t)=\eta'_i(|z|)\tau'_i(t)$ as in Definition
1. We further assume that $\varrho_0=\varrho$, $\varrho'_0=\varrho'$,  
$\varrho_{i+1}=1$ in a neighborhood of the support of
$\varrho_i$, and $\varrho'_{i+1}=1$ in a neighborhood of the support of  
$\varrho'_i$ and that the $\eta_i(|z|)$ are monotone
decreasing in $|z|$.  We also choose
$\{\gamma^+_i\}$ and $\{\gamma^0_i\}$ such that $\gamma_i^+\in\mathcal  
G^+$, $\gamma^+_{i+1}=1$, and $\gamma_i^0\in\mathcal
G^0$ and $\gamma^0_{i+1}=1$ in neighborhoods of the supports of  
$\gamma^+_i$ and $\gamma^0_i$, respectively. Further we
require that $\gamma^0_i=1$ in a neighborhood of the support of  
derivatives of $\gamma^+_i$. Substituting
$\Psi^{-s_1}\Gamma^+_0u$  for $u$ in (\ref{hypest}), replacing  
$s_0+s_1$ by $s_0$, we  have
$$
\|\varrho\Psi^{-s_1}\Gamma^+_0u\|_{s_1}\le C(\|\varrho'  
E_k\Psi^{-s_1}\Gamma^+_0u\|_{s_2}+\|u\|_{-s_0}).
$$
Since $\gamma^+_1\varrho\gamma^+_0=\varrho\gamma^+_0$, 
\begin{align*}
\|\varrho\Psi^{- 
s_1}\Gamma^+_0u\|_{s_1}&=\|\Psi^{s_1}\Gamma^+_1\varrho\Psi^{- 
s_1}\Gamma^+_0u\|+
O(\|u\|_{-1})\\
&=\|\varrho\Gamma^+_0u\|+\|\Psi^{s_1}\Gamma^+_1[\varrho,\Psi^{s_1}]
\Gamma^+_0u\|+O(\|u\|_{-s_0}).
\end{align*}
Furthermore, $\Psi^{s_1}\Gamma^+_1[\varrho,\Psi^{s_1}]\Gamma^+_0$ is an  
operator of order $-1$; hence we get
$$
\|\Psi^{s_1}\Gamma^+_1[\varrho,\Gamma^+_0]\Psi^{-s_1}u\|\le  
C(\|u\|_{-s_0})
$$
and
$$
\|\varrho\Gamma^+_0u\|\le
C(\|\varrho'E_k\Psi^{-s_1}\Gamma^+_0u\|_{s_2}+\|u\|_{-s_0}).
$$
Next we have
$$
\|\varrho'  
E_k\Psi^{-s_1}\Gamma^+_0u\|_{s_2}\le\|\Psi^{s_2- 
s_1}\Gamma^+_0\varrho'E_ku\|+
\|[\varrho'E_k,\Psi^{-s_1}\Gamma^+_0]u\|_{s_2}.
$$
Since the symbol of $\gamma^0_1\gamma^+_1\varrho'_1=1$ in a  
neighborhood of the symbol of\break
$[\varrho'E_k,\Psi^{-s_1}\Gamma^+_0]$ and since the order of  
$[\varrho'E_k,\Psi^{-s_1}\Gamma^+_0]$ is $-s_1+1$, we have
$$
\|[\varrho'E_k,\Psi^{-s_1}\Gamma^+_0]u\|_{s_2}\le  
C(\|\varrho'_1\Gamma^0_1\Gamma^+_1u\|_{s_2-s_1+1}+\|u\|_{-s_0}).
$$
Applying Proposition 3, we have
$$
\|\varrho'_1\Gamma^0_1\Gamma^+_1u\|_{s_2-s_1+1}\le  
C(\|\varrho'_1E_k\Gamma^+_1u\|_{s_2-s_1-1}+
\|\varrho'_2\Gamma^0_2\Gamma^+_1u\|_{s_2-s_1}+\|u\|_{-\infty})
$$
so that
\begin{eqnarray*}
&&\hskip-6pt\|\varrho' E_k\Psi^{-s_1}\Gamma^+_0u\|_{s_2}\\
&&\quad\le  
C(\|\Psi^{s_2-s_1}\Gamma^+_0\varrho'E_ku\|+
\|\varrho'_1\Gamma^0_2\Gamma^+_1u\|_{s_2- 
s_1}+\|\varrho'_1E_k\Gamma^+_1u\|_{s_2-s_1-1}+\|u\|_{-s_0}).
\end{eqnarray*}
Therefore
\begin{eqnarray*} 
\|\varrho\Gamma^+_0u\| &\le& C(\|\Psi^{s_2-s_1}\Gamma^+_0\varrho'E_ku\|\\
&&\phantom{ C(} +
\|\varrho'_1E_k\Gamma^+_1u\|_{s_2-s_1- 
1}+\|\varrho'_1E_k\Gamma^+_1u\|_{s_2-s_1-1}+\|u\|_{-s_0}).
\end{eqnarray*}
Now we have
$$
\|\varrho'_1E_k\Gamma^+_1u\|_{s_2-s_1-1}\le\|\Psi^{s_2-s_1- 
1}\Gamma^+_1\varrho'_1E_ku\|+
\|[\varrho'_1E_k,\Gamma^+_1]u\|_{s_2-s_1-1},
$$
again since $[\varrho'_1E_k,\Gamma^+_1]$ is an operator of order one  
and since $\varrho'_22\gamma^0_2\gamma^+_2=1$ in a
neighborhood of its symbol, we get
$$
\|[\varrho'_2E_k,\Gamma^+_1]u\|_{s_2-s_1-1}\le  
C(\|\varrho'_2\Gamma^0_2\Gamma^+_2u\|_{s_2-s_1}+\|u\|_{-s_0}).
$$
Then, again applying Proposition 3, we have
$$
\|\varrho'_2\Gamma^0_2\Gamma^+_2u\|_{s_2-s_1}\le  
C(\|\varrho'_3\Gamma^0_3E_k\Gamma^+_2u\|_{s_2-s_1-2}+\|u\|_{-\infty})
$$
so that
$$
\|\varrho'_2E_k\Gamma^+_1u\|_{s_2-s_1-1}\le  
C(\|\Psi^{s_2-s_1-1}\Gamma^+_1\varrho'_2E_ku\|+
\|\varrho'_3\Gamma^0_3E_k\Gamma^+_2u\|_{s_2-s_1-2}+\|u\|_{-s_0}).
$$
Hence
\begin{eqnarray*}
\|\varrho\Gamma^+u\|&\le& C(\|\Psi^{s_2-s_1}\Gamma^+_0\varrho'E_ku\|\\
&&\phantom{C(} +
\|\Psi^{s_2-s_1- 
1}\Gamma^+_1\varrho'_2E_ku\|+\|\varrho'_4\Gamma^0_3E_k\Gamma^+_2u\|_{s_2 
-s_1-2}+\|u\|_{-s_0}).
\end{eqnarray*}
Proceeding inductively we obtain
\begin{eqnarray*}
\|\varrho\Gamma^+u\|&\le  &
C\Big(\sum_{i=0}^N\|\Psi^{s_2-s_1-i}\Gamma^+_i\varrho_i'E_ku\|\\
&&\phantom{C\Big(}+
\|\varrho'_{N+3}\Gamma^0_{N+2}E_k\Gamma^+_{N+1}u\|_{s_2-s_1-N- 
1}+\|u\|_{-s_0}\Big).
\end{eqnarray*}
Since $\|[\Psi^{s_2-s_1-i}\Gamma^+_i,\eta_i']E_ku\|$ can be  
incorporated in the successive terms, we get, by choosing
$N\ge s_2-s_1+1-s_0$
$$
\|\varrho\Gamma^+u\|\le  
C\Big(\sum_{i=0}^N\|\Psi^{s_2-s_1-i}\Gamma^+_i\tau'_iE_ku\|+\|u\|_{-s_0}\Big).
$$
Let $\tilde\tau\in C^\infty_0$ with $\tilde\tau=1$ on the support of  
$\tau'_N$; then $\tau'_iE_ku=\tau'_iE_k\tilde\tau u$
when $i\le N$ so that replacing $u$ by $\tilde\tau u$ we obtain
$$
\|\varrho\Gamma^+\tilde\tau u\|\le  
C\Big(\sum_{i=0}^N\|\Psi^{s_2-s_1-i}\Gamma^+_i\tau'_iE_ku\|+
\|\tilde\tau u\|_{-s_o}\Big).
$$
Hence, since $\gamma^0=1$ in a neighborhood of the support of
the symbol of $[\Gamma^+,\tilde\tau]$ and thus can be incorporated in  
the estimate as above,   we have
\begin{align*}
\|\varrho\Gamma^+u\|^2\le  
C\Big(\sum_{i=0}^N\|\Psi^{s_2-s_1-i}\Gamma^+_i\tau'_iE_ku\|^2+\|\tilde\tau  
u\|^2_{-s_o}\Big).
\end{align*}
Choosing $\tilde\gamma^+$ so that $\tilde\gamma^+=1$ in a neighborhood  
of the supports of the $\gamma_i^+$, we have
$\tilde\tau\tilde\gamma^+=1$ in a neighborhood of the support of the  
symbol of $\Psi^{s_2-s_1-i}\Gamma^+_i\tau'_iE_k$.
Then we obtain
\begin{align*}
\|\varrho\Gamma^+u\|^2&\le  
C\Big(\sum_{i=0}^N\|\Psi^{s_2-s_1-i}\tilde\Gamma^+\tilde\tau  
E_ku\|^2+\|\tilde\tau u\|^2_{-s_o}\Big)\\
&\le C(\|\Psi^{s_2-s_1}\tilde\Gamma^+\tilde\tau E_ku\|^2+\|\tilde\tau  
u\|^2_{-s_o}).
\end{align*}
We define $h_{\lambda}$ by
$$
h_{\lambda}(z,t)={\rm exp}\left(-\lambda^2(|z|^2-it)\right),
$$
since $\tilde\tau h_{\lambda}\in\mathcal S$ and obtain
$$
\|\varrho\Gamma^+h_\lambda\|\le  
C\left(\|\Psi^{s_2-s_1}\tilde\Gamma^+_i\tilde\tau E_kh_{\lambda}\|+
\|\tilde\tau h_{\lambda}\|_{-s_o})\right).
$$
Assuming that $\eta(|z|)$ is monotone decreasing we have  
$\eta(|z|)\ge\eta(\lambda|z|)$; hence, setting
$\eta_\lambda(z)=\eta(\lambda|z|)$, we obtain

$$
\|\varrho\Gamma^+h_{\lambda}\|\ge\|\eta_\lambda\tau\Gamma^+h_{\lambda}\|.
$$
Then, setting $x'=(x_1,x_2)$, $y'=(y_1,y_2)$, and $\xi'=(\xi_1,\xi_2)$  
and changing variables $\lambda y'\to y'$,
$\xi'\to\lambda\xi'$, and $\xi'_3\to\xi'_3+\lambda^2$, we get
\begin{align*}
&\eta_\lambda\tau\Gamma^+h_{\lambda}(x)=
\int\exp(i(x-y)\cdot\xi)\tau(y_3)\gamma^+(\xi)\exp(-\lambda^2(|y'|^2- 
iy_3))dyd\xi\\
&=\int\exp(i(x'-y')\cdot\xi'+x_3\xi_3-y_3(\xi_3- 
\lambda^2))\tau(y_3)\gamma^+(\xi)\exp(-\lambda^2|y'|^2)dyd\xi\\
&=\lambda^{-2}\int\exp(i(\lambda  
x'-y')\cdot\xi')\\
&\quad+(x_3- 
y_3)\xi_3)\tau(y_3)\gamma^+(\lambda\xi',\lambda^2+\xi_3)
\exp(-|y'|^2)dyd\xi.
\end{align*}
Making the change of variables $\lambda x'\to x'$ we have
\begin{align*}
&\|\eta_\lambda\tau\Gamma^+h_{\lambda,\delta}\|^2\\
&\quad=\frac{1}{\lambda^6}\int 
|
\int\exp(i(x- 
y)\cdot\xi)\tau(y_3)\gamma^+(\lambda\xi',\lambda^2+\xi_3)\exp(- 
|y'|^2)dyd\xi|^2dx.
\end{align*}
Given $(\xi',\xi_3)\in {\rm supp}(\gamma^+)$ we have
$$
\lim_{\lambda\to\infty}|\frac{\lambda\xi'}{\xi_3+\lambda^2}|=0,
$$
and   there exists $\tilde\lambda$ such that
$\gamma^+(\lambda\xi',\lambda^2+\xi_3)=1$ when  
$\lambda\ge\tilde\lambda$.
Hence we have
$\lim_{\lambda\to\infty}\gamma^+(\lambda\xi',\lambda^2+\xi_3)=1$; thus  
there  exist $\lambda_0$ such that
$$
\|\eta_\lambda\tau\Gamma^+h_{\lambda}\|^2\ge\frac{1}{2\lambda^6}
\int|\int\exp(i(x-y)\cdot\xi)\tau(y_3)\exp(-|y'|^2)dyd\xi|^2dx,
$$
when $\lambda\ge\lambda_0$, therefore there exists $C$ independent of  
$\lambda$ such that
$$
\|\varrho\Gamma^+h_{\lambda}\|\ge\frac{C}{\lambda^3},
$$
when $\lambda\ge\lambda_0$.
 
 Next, we will estimate the term $\|\tilde\tau  
h_{\lambda}\|_{-s_o}$. We will use the facts that
$\frac{1}{m!}\bar L^m(\bar z^m)\break=1$ and that $\bar L(h_\lambda)=0$.  
Taking $m\le s_o$, we
have
\begin{align*}
&\mathcal F(\Lambda^{-s_o}\tilde\tau h_{\lambda})(\xi)=
\int(1+|\xi|^2)^{\frac{-s_o}{2}}\tilde\tau(x_3)\exp(-i(x\cdot\xi- 
\lambda^2x_3)-\lambda^2|z|^2)dx\\
&\quad=\frac{1}{m!}\int\bar L^m(\bar z^m)(1+|\xi|^2)^{\frac{-s_o}{2}}
\tilde\tau)\exp(-i(x\cdot\xi-\lambda^2x_3)-\lambda^2|z|^2)dx\\
&\quad=-\frac{1}{m!}\int\bar z^m(1+|\xi|^2)^{\frac{-s_o}{2}}
\bar  
L^m(\tilde\tau(x_3)\exp(-i(x\cdot\xi-\lambda^2x_3)-\lambda^2|z|^2)dx\\
&\quad=-\frac{1}{m!}\int\bar z^m(1+|\xi|^2)^{\frac{-s_o}{2}}
\exp(-\lambda^2(|z|^2-ix_3))\bar  
L^m(\tilde\tau(x_3)\exp(-ix\cdot\xi))dx\\
&\quad=-\frac{1}{m!}\int\bar z^m(1+|\xi|^2)^{\frac{-s_o}{2}}
\exp(-\lambda^2(|z|^2))\bar  
L^m(\tilde\tau(x_3)\\
&\hskip.4in\cdot\exp(-ix'\cdot\xi'-ix_3(\xi_3-\lambda^2))dx
\end{align*}
and
\begin{align*}
&\bar  
L^m\left(\tilde\tau(x_3)\exp(-ix'\cdot\xi'+i-ix_3(\xi_3- 
\lambda^2))\right)\\
&=\sum_{j=0}^ma_j(x_3)z^j(i\xi_1+\xi_2-2z\xi_3)^{m-j}\exp(- 
ix'\cdot\xi'-ix_3(\xi_3-\lambda^2)).
\end{align*}
Thus, setting  
$w^{(m)}(x,\xi)=\sum_{j=0}^ma_j(x_3)z^j(i\xi_1+\xi_2-2z\xi_3)^{m-j}$  
and denoting the corresponding
pseudodifferential operator by $W^{(m)}$, we have
$$
\|\tilde\tau h_{\lambda}\|_{-s_o}=C\|W^{(m)}\bar  
z^mh_{\lambda}\|_{-s_o}\le C\|z^m\tilde\tau'h_{\lambda}\|_{m-s_o}
\le C\|z^m\tilde\tau'h_{\lambda}\|,
$$
where $\tilde\tau'\in C^\infty_0(\mathbb R)$ and $\tilde\tau'=1$ in a  
neighborhood of the support of $\tilde\tau$. Now, changing
coordinates $\lambda z\to z$, we get
$$
\|z^m\tilde\tau'h_{\lambda}\|^2=\int|z|^{2m}\tilde\tau'(x_3)^2\exp(- 
2\lambda^2|z|^2)dx\le\frac{C}{\lambda^{2m+2}}.
$$
To estimate the remaining terms we have
$$
E_kh_\lambda(z,t)=-2(k+1)\lambda^2|z|^{2k}h_\lambda(z,t).
$$
Therefore, with the coordinate change $\lambda x'\to x'$, we get
\begin{align*}
&
\mathcal F\left(\Psi^{s}\Gamma^+\tau E_kh_\lambda\right)(\xi)\\
&\qquad=
C\mathcal  
F\left(\Psi^{s}\Gamma^+\tau\lambda^2|z|^{2k}h_\lambda\right)(\xi)\\
&\qquad=C\lambda^{-2}(1+\xi_3^2)^{\frac{s}{2}}\gamma^+(\xi)\mathcal
F\left(\tau(x_3)|z|^{2k}\exp(-\lambda^2|z|^2)\right)(\xi',\xi_3- 
\lambda^2)\\
&\qquad=C\lambda^{-2k-2}(1+\xi_3^2)^{\frac{s}{2}}\gamma^+(\xi)\hat\tau(\xi_3- 
\lambda^2)\mathcal
F\left(|z|^{2k}\exp(-2|z|^2)\right)(\lambda^{-1}\xi').
\end{align*}
Then, integrating and making the changes of coordinates  
$\xi'\to\lambda\xi'$,  $\xi_3\to\xi_3+\lambda^2$, we get
\begin{align*}
&\hskip-9pt\|\Psi^{s}\Gamma^+\tau E_kh_\lambda\|^2\\
&\le  
C\lambda^{-4k-4}\int(1+\xi_3^2)^{s}
\gamma^+(\xi)\hat\tau(\xi_3-\lambda^2)|^2|\mathcal  
F\left(|z|^{2k}\exp(-|z|^2)\right)(\lambda^{-1}\xi')|^2d\xi\\
&\le  
C\lambda^{-4k- 
2}\\&\enspace
\cdot\int(1+(\xi_3+\lambda^2)^2)^{s}\gamma^+(\lambda\xi',\xi_3+\lambda^2) 
\hat\tau(\xi_3)|^2|\mathcal
F\left(|z|^{2k}\exp(-|z|^2)\right)(\xi')|^2d\xi.
\end{align*}
Then if $s\ge0$ and if $\lambda$ is sufficiently large we have
$$
\|\Psi^{s}\Gamma^+\tau E_kh_\lambda\|^2\le C\lambda^{4s-4k-2}.
$$
We assume $k\ge 1$; if  $s_2-s_1<0$ then  
$$
\|\Psi^{s_1-s_2}\varrho\Gamma^+h_\lambda\|^2\le  
C\left(\|\tilde\Gamma^+_i\tilde\tau E_kh_{\lambda}\|^2+
\|\tilde\tau h_{\lambda}\|^2_{-s_o})\right)\le C\lambda^{-4k-2}
$$
and, by Lemma 1,  
\begin{align*}
\|\Psi^{s_1-s_2}\varrho\Gamma^+h_\lambda\|^2&=
\|\Psi^{s_1- 
s_2}\eta_\lambda\tau\Gamma^+h_{\lambda}\|^2+O(\|\tau\Gamma^+h_{\lambda}\ 
|^2_{-s_o})\\
&\ge  
C\lambda^{2s_2-2s_1}\|\eta_\lambda\tau\Gamma^+h_{\lambda}\|^2- 
C'(\|\tau\Gamma^+h_{\lambda}\|^2_{-s_o})\\
&\ge C(\lambda^{2s_2-2s_1-2}+\lambda^{-2m-2}).
\end{align*}
This implies that for large $\lambda$ we have $\lambda^{2s_2-2s_1-2}\le  
C(\lambda^{-4k-2}+\lambda^{-2m-2})$,
which is a contradiction, so that $s_2-s_1\ge0$ and    
\begin{align*}
C_1\lambda^{-6}&\le C_2\|\varrho\Gamma^+h_\lambda\|^2\le  
C\left(\|\Psi^{s_2-s_1}\tilde\Gamma^+_i\tilde\tau E_kh_{\lambda}\|^2+
\|\tilde\tau h_{\lambda}\|^2_{-s_o}\right)\\*
&\le C_3(\lambda^{4s_2-4s_1-4k-2}+\lambda^{-2m-2}).
\end{align*}
Therefore, if $m$ large we get $C_1\le  
2C_3\lambda^{4(s_2-s_1-k+1)}$ for large $\lambda$. Hence
$s_2-s_1-k+1\ge0$, which concludes the proof of the proposition and  
also of Theorem~B.

\section{Elliptic and subelliptic microlocalizations}

In this section we will show that the {\it a priori\/} estimates for the  
operator $E_k$ gain two derivatives in the $0$ microlocalization
and gains one derivative in the $-$ microlocalization, these gains are  
in the Sobolev norms. Without loss of generality we will
deal only with microlocalizations near the origin, taking  
$\alpha=0$ and setting $\mathfrak G^0=\mathfrak G^0_0$ and
$\mathfrak G^-=\mathfrak G^-_0$. The subscript $\alpha$ will be  
dropped from the corresponding operators.
 
\begin{prop} Let $U$ and $U'$ be neighborhoods of the origin with $\bar  
U\subset U'$ and $|z|\le a$ on $U'${\rm ,} where $a$ is
sufficiently small as in Lemma {\rm 3}. Suppose that $\varrho\in  
C_0^\infty(U)$ and $\varrho'\in C_0^\infty(U')$ with
$\varrho'=1$ on a neighborhood of $\bar U$. Further suppose that  
$\gamma^0,\tilde\gamma^0\in\mathcal G^0$ with $\tilde\gamma^0=1$
on a neighborhood of the support of $\gamma^0$. Then{\rm ,} given  
$s,s_0\in\mathbb R${\rm ,} there exists
$C=C(\varrho,\varrho',\gamma^0,\tilde\gamma^0,s,s_0)$ such that
$$
\|\varrho\Gamma^0u\|^2_{s+2}+\|\varrho\Gamma^0\bar  
z^kLu\|^2_{s+1}+\|\varrho\Gamma^0\bar Lu\|^2_{s+1}\le
C(\|\varrho'\tilde\Gamma^0E_ku\|^2_s+\|u\|^2_{-s_0}),
$$
for all $u\in\mathcal S${\rm ,} where $\mathcal S$ denotes the Schwartz class  
of rapidly decreasing functions.
\end{prop}

  \Proof  Let $\{\varrho_i\}$ be a sequence of functions  
such that
$\varrho_i\in C_0^\infty(U)$, $\varrho_0=\varrho$, $\varrho_{i+1}=1$ in  
a neighborhood of
the support of $\varrho_i$, and such that $\varrho'=1$ in a  
neighborhood of the supports of all the $\varrho_i$. Let
$\{\gamma^0_i\}$ be a sequence in $\mathcal G^0$ such that  
$\gamma^0_0=\gamma^0$, $\gamma^0_{i+1}=1$ in a neighborhood of
the support of $\gamma_i$, and $\tilde\gamma^0=1$ in a neighborhood of  
the supports of all the $\gamma^0_i$. Then
substituting $\varrho\Lambda^{s+1}\Gamma^0_1u$ for $u$ in Lemma 3 we  
have
$$
\|\Gamma^0\varrho\Lambda^{s+1}\Gamma^0_1u\|^2_1\le C(\|\Gamma^0\bar  
L\varrho\Lambda^{s+1}\Gamma^0_1u\|^2+
\|\varrho\Lambda^{s+1}\Gamma^0_1u\|^2).
$$
Hence
$$
\|\Gamma^0\varrho\Lambda^{s+1}\Gamma^0_1u\|^2_1\le C(\|\Gamma^0\bar  
z^kL\varrho\Lambda^{s+1}\Gamma^0_1u\|^2+
\|\Gamma^0\bar  
L\varrho\Lambda^{s+1}\Gamma^0_1u\|^2+\|\varrho\Lambda^{s+1}\Gamma^0_1u\| 
^2).
$$
Then
\begin{eqnarray*}
 \|\Gamma^0\varrho\Lambda^{s+1}\Gamma^0_1u\|^2_1&=&\|\varrho\Gamma^0u\|_{s 
+2}^2+
O(\|\Lambda^1[\Gamma^0\varrho,\Lambda^{s+1}]\Gamma^0_1u\|^2+\|\Lambda^{s 
+2}[\Gamma^0,\varrho]\Gamma^0_1u\|^2
\\
&&\phantom{\|\varrho\Gamma^0u\|_{s 
+2}^2+
O(}+\|\Lambda^{s+2}\varrho(\Gamma^0\Gamma^0_1-\Gamma^0)u\|^2).
\end{eqnarray*}
Since $[\Gamma^0\varrho,\Lambda^{s+1}]$ is a pseudodifferential  
operator of order $s+1$ and since $\varrho_1=1$ on the
support of its symbol, we have
$$
\|\Lambda^1[\Gamma^0\varrho,\Lambda^{s+1}]\Gamma^0_1u\|^2\le  
C(\|\varrho_1\Gamma^0_1u\|^2_{s+1}+\|\Gamma^0_1u\|^2_{-\infty}).
$$
The operator	$[\Gamma^0,\varrho]$ is of order $-1$ and $\varrho_1=1$ on  
the support of its symbol, so that
$$
\|\Lambda^{s+2}[\Gamma^0,\varrho]\Gamma^0_1u\|^2\le  
C(\|\varrho_1\Gamma^0_1u\|^2_{s+1}+\|\Gamma^0_1u\|^2_{-\infty}).
$$
The symbol of the operator  
$\Lambda^{s+2}\varrho(\Gamma^0\Gamma^0_1-\Gamma^0)$ is zero so that
$$
\|\Lambda^{s+2}\varrho(\Gamma^0\Gamma^0_1-\Gamma^0)u\|^2\le  
C\|u\|^2_{-\infty}.
$$
Then we obtain
$$
\|\varrho\Gamma^0u\|_{s+2}^2\le
C(\|\Gamma^0\varrho\Lambda^{s+1}\Gamma^0_1u\|^2_1+\|\varrho_1\Gamma^0_1u 
\|^2_{s+1}+\|u\|^2_{-\infty}),
$$
so that
$$
\|\varrho\Gamma^0u\|_{s+2}^2\le C(\|\Gamma^0\bar  
z^kL\varrho\Lambda^{s+1}\Gamma^0_1u\|^2+
\|\Gamma^0\bar  
L\varrho\Lambda^{s+1}\Gamma^0_1u\|^2+\|\varrho_1\Gamma^0_1u\|^2_{s+1}+\| 
u\|^2_{-\infty}).
$$
The following lemma which involves a vector field $X$ will be applied  
with $X=\bar z^kL$ and $X=\bar L$.

\begin{lemma}If $X$ is a complex vector field on $\mathbb R^3$ then
$$
\|\Gamma^0X\varrho\Lambda^{s+1}\Gamma^0_1u\|^2\le  
C(\|\varrho\Gamma^0Xu\|^2_{s+1}
+\|\varrho_1\Gamma^0_1u\|^2_{s+1}+\|u\|^2_{-\infty})
$$
and
\begin{eqnarray*}
\|\varrho\Gamma^0Xu\|^2_{s+1}&=&(\Lambda^s\varrho\Gamma^0X^*Xu,\Lambda^{s 
+2}\varrho\Gamma^0u)\\
&&+
O(\|\varrho_1\Gamma^0_1u\|^2_{s+1}+\|\varrho\Gamma^0u\|^2_{s+2}\|\varrho 
_1\Gamma^0_1u\|^2_{s+1}\\
&&\phantom{+O(}+\|\varrho_1\Gamma^0_1Xu\|^2_s+\|\varrho_1\Gamma^0_1u\|_{s+1}\| 
\varrho\Gamma^0Xu\|_{s+1}+\|u\|^2_{-\infty}),
\end{eqnarray*}
for all $u\in\mathcal S$.
\end{lemma}

\Proof  We have
$$
\|\Gamma^0X\varrho\Lambda^{s+1}\Gamma^0_1u\|\le\|\Gamma^0\varrho\Lambda^ 
{s+1}\Gamma^0_1Xu\|+
\|\Gamma^0[X,\varrho\Lambda^{s+1}\Gamma^0_1]u\|.
$$
The operator  
$P=\Gamma^0\varrho\Lambda^{s+1}\Gamma^0_1X- 
\Lambda^{s+1}\Gamma^0_1\varrho\Gamma^0X$ is of order
$s+1$ and $\varrho_1\gamma^0_1=1$ in a neighborhood of the symbol of  
$P$; hence
$$
\|Pu\|\le C(\|P\varrho_1\Gamma^0_1u\|+\|u\|_{-\infty})\le  
C(\|\varrho_1\Gamma^0_1u\|_{s+1}+\|u\|_{-\infty}).
$$
Since $\gamma^0_1=1$ in a neighborhood of the support of the symbol of  
$\varrho\Gamma^0$, we get
$$
\|\Gamma^0\varrho\Lambda^{s+1}\Gamma^0_1Xu\|\le  
C(\|\varrho\Gamma^0Xu\|_{s+1}+\|u\|_{-\infty}).
$$
Furthermore, $\Gamma^0[X,\varrho\Lambda^{s+1}\Gamma^0_1]$ is of order  
$s+1$ and $\varrho_1\gamma^0_1=1$ in a neighborhood of
the support of its symbol so that
$$
\|\Gamma^0[X,\varrho\Lambda^{s+1}\Gamma^0_1]u\|\le  
C(\|\varrho_1\Gamma^0_1u\|_{s+1}+\|u\|_{-\infty}),
$$
which proves the first part of the lemma.

  For the second part of the lemma we write
\begin{align*}
\|\varrho\Gamma^0Xu\|^2_{s+1}&=(\Lambda^{s+1}\varrho\Gamma^0Xu,\Lambda^{ 
s+1}\varrho\Gamma^0Xu)\\
&=(\Lambda^{s+1}\varrho\Gamma^0Xu,[\Lambda^{s+1}\varrho\Gamma^0,X]u)
+([X^*,\Lambda^{s+1}\varrho\Gamma^0]Xu,\Lambda^{s+1}\varrho\Gamma^0u)\\
&\quad+(\Lambda^s\varrho\Gamma^0X^*Xu,\Lambda^{s+2}\varrho\Gamma^0u).
\end{align*}
Then, since $[\Lambda^{s+1}\varrho\Gamma^0,X]$ is of order $s+1$ and  
$\varrho_1\gamma^0_1=1$ in a neighborhood of its symbol,
$$
\|[\Lambda^{s+1}\varrho\Gamma^0,X]u\|^2\le  
C(\|\varrho_1\Gamma^0_1u\|_{s+1}^2+\|u\|_{-\infty}).
$$
Then
\begin{eqnarray*}
&&\hskip-7pt
([X^*,\Lambda^{s+1}\varrho\Gamma^0]Xu,\Lambda^{s+1}\varrho\Gamma^0u)\\
&&=
([(\Lambda^{s+1}\varrho\Gamma^0)^*,[X^*,\Lambda^{s+1}\varrho\Gamma^0]]Xu 
,u) +((\Lambda^{s+1}\varrho\Gamma^0)^*Xu,[X^*,\Lambda^{s+1}\varrho 
\Gamma^0]^*u).
\end{eqnarray*}
Let  
$Q=[(\Lambda^{s+1}\varrho\Gamma^0)^*,[X^*,\Lambda^{s+1}\varrho\Gamma^0]] 
$; then $Q$ has order $2s+1$ and
$\varrho_1\gamma^0_1=1$ in a neighborhood of its symbol. Thus
\begin{align*}
|(QXu,u)|&\le  
C(|(Q\varrho_1\Gamma^0_1Xu,\varrho_1\Gamma^0_1u)|+\|u\|_{-\infty}^2)\\
&\le  
C(\|\varrho_1\Gamma^0_1Xu\|^2_{s}+\|\varrho_1\Gamma^0_1u\|^2_{s+1}+\|u\| 
_{-\infty}^2).
\end{align*}
The symbol of the operator  
$(\Lambda^{s+1}\varrho\Gamma^0)^*-\Lambda^{s+1}\varrho\Gamma^0$ is  
zero, the order of\break
$[X^*,\Lambda^{s+1}\varrho\Gamma^0]^*$ is $s+1$ and  
$\varrho_1\gamma^0_1=1$ on a neighborhood of its support. Hence
$$
|((\Lambda^{s+1}\varrho\Gamma^0)^*Xu,[X^*,\Lambda^{s+1}\varrho\Gamma^0]^ 
*u)|\le
C(\|\varrho\Gamma^0Xu\|_{s+2}\|\varrho_1\Gamma_1^0u\|_{s+1}+\|u\|_{- 
\infty}^2).
$$
Combining these we conclude the proof of the lemma.

 Returning to the proof of the proposition, by using the above  
lemma, when $X=\bar z^kL$ and when
$X=\bar L$, we obtain
\begin{multline*}
 \|\varrho\Gamma^0u\|^2_{s+2}+\|\varrho\Gamma^0\bar  
z^kLu\|^2_{s+1}+\|\varrho\Gamma^0\bar Lu\|^2_{s+1}\\
\le
C(\|\varrho\Gamma^0 E_ku\|^2_s+\|\varrho_1\Gamma_1^0\bar z^kLu\|^2_{s} 
 +\|\varrho_1\Gamma_1^0\bar Lu\|^2_{s}+\|u\|_{-\infty}^2).
\end{multline*}
Replacing $\varrho$ by $\varrho_i$, $\varrho_1$ by $\varrho_{i+1}$,  
$\Gamma^0$ by $\Gamma^0_i$,
$\Gamma_1^0$ by $\Gamma_{i+1}^0$, and $s$ by $s-i$ we obtain
\begin{multline*}
 \|\varrho_i\Gamma_i^0u\|^2_{s+2-i}+\|\varrho_i\Gamma_i^0\bar  
z^kLu\|^2_{s+1-i}+\|\varrho_i\Gamma_i^0\bar Lu\|^2_{s+1-i}\\
\le
C(\|\varrho_i\Gamma_i E_ku\|^2_{s-i}+\|\varrho_{i+1}\Gamma_{i+1}^0\bar  
z^kLu\|^2_{s-i} +\|\varrho^0_{i+1}\Gamma_{i+1}^0\bar  
Lu\|^2_{s-i}+\|u\|_{-\infty}^2).
\end{multline*}
Proceeding inductively, we obtain
\begin{multline*}
 \|\varrho\Gamma^0u\|^2_{s+2}+\|\varrho\Gamma^0\bar  
z^kLu\|^2_{s+1}+\|\varrho\Gamma^0\bar Lu\|^2_{s+1}\\
\le
C\Big(\sum_{i=0}^N\|\varrho_i\Gamma^0_i E_ku\|^2_{s-i} 
+\|\varrho_{N+1}\Gamma_{N+1}^0\bar  
z^kLu\|^2_{s-N}+\|\varrho^0_{i+N}\Gamma_{i+N}^0\bar Lu\|^2_{s-N}+
\|u\|_{-\infty}^2\Big).
\end{multline*}
Setting $N\ge s_0+s+1$ we conclude the proof of the proposition since  
$$\|\varrho_i\Gamma^0_i E_ku\|^2_{s-i}\le
C(\|\varrho'\tilde\Gamma^0 E_ku\|^2_{s}+\|u\|^2_{-\infty}).
$$

\begin{prop}\hskip-2pt Given neighborhoods of the origin $U$ and $U'$ with $\bar  
U\!\subset\! U'${\rm ;} suppose that $\varrho\in C_0^\infty(U)$
and $\varrho'\in C_0^\infty(U')$ with $\varrho'=1$ on a neighborhood of~$\bar U$. Further
suppose that
$\gamma^-,\tilde\gamma^-\in\mathcal G^-$ with $\tilde\gamma^-=1$ on a  
neighborhood of the support of $\gamma^-$. Then{\rm ,} given
$s,s_0\in\mathbb R${\rm ,} there exists
$C=C(\varrho,\varrho',\gamma^-,\tilde\gamma^-,s,s_0)$ such that
$$
\|\varrho\Gamma^-u\|^2_{s+1}+\|\varrho\Gamma^-\bar  
z^kLu\|^2_{s+\frac{1}{2}}+
\|\varrho\Gamma^-\bar Lu\|^2_{s+\frac{1}{2}}\le  
C(\|\varrho'\tilde\Gamma^-E_ku\|^2_s+\|u\|^2_{-s_0}),
$$
for all $u\in\mathcal S$.
\end{prop}

  \Proof  The proof is entirely analogous to that of the  
above proposition.  We use Lemma 5 in place of Lemma 3
and substitute $\varrho\Lambda^{s+\frac{1}{2}}\Gamma^-_1u$ for $u$ we  
obtain
$$
\|\Gamma^-\varrho\Lambda^{s+\frac{1}{2}}\Gamma^- 
_1u\|^2_{\frac{1}{2}}\le C(\|\Gamma^0\bar
L\varrho\Lambda^{s+\frac{1}{2}}\Gamma^-_1u\|^2+
\|\varrho\Lambda^{s+\frac{1}{2}}\Gamma^0_1u\|^2).
$$
Then one proceeds exactly as above to obtain the proof.

   In the case $k=0$ the vectorfields $L$ and $\bar L$ play  
exactly the same role and so we obtain the following.
\begin{prop} Given neighborhoods of the origin $U$ and $U'$ with $\bar  
U\!\subset\! U'$. Suppose that $\varrho\in C_0^\infty(U)$
and $\varrho'\in C_0^\infty(U')$ with $\varrho'=1$ on a neighborhood of~$\bar U$. Further
suppose that
$\gamma^+,\tilde\gamma^+\in\mathcal G^+$ with $\tilde\gamma^+=1$ on a  
neighborhood of the support of $\gamma^+$. Then{\rm ,} given
$s,s_0\in\mathbb R${\rm ,} there exists
$C=C(\varrho,\varrho',\gamma^+,\tilde\gamma^+,s,s_0)$ such that
$$
\|\varrho\Gamma^+u\|^2_{s+1}+\|\varrho\Gamma^+Lu\|^2_{s+\frac{1}{2}}+
\|\varrho\Gamma^+\bar Lu\|^2_{s+\frac{1}{2}}\le  
C(\|\varrho'\tilde\Gamma^+E_0u\|^2_s+\|u\|^2_{-s_0}),
$$
for all $u\in\mathcal S$.
\end{prop}

\section{The operator $E_0$ and gain of derivatives}

Since $E_0$ is a real operator, it can be written as $E_0=-X^2-Y^2$, where  
$X=\frac{1}{\sqrt2}\Re L$ and $Y=\frac{1}{\sqrt2}\Im L$.
Thus it is one of the simplest operators that\break\vskip-12pt\noindent satisfy H\"ormander's  
condition and it is well understood. Nevertheless, it is
instructive to write it in terms of $L$ and $\bar L$ and analyze it  
microlocally in the framework of the previous section. The
operator $E_0$ gains one derivative. As we have seen the operators  
$E_k$ do not gain derivatives when $k>0$ and $z=0$; in a
neighborhood on which $z\neq0$ they do gain derivatives and they also  
gain in the $0$ and $-$ microlocalizations.

  In the analysis of $E_0$ we can assume, without loss of  
generality, that $\alpha=0$ and we set
$\gamma=\gamma_0$, and $\Gamma=\Gamma_0$. The  basic observation is  
that the gain of derivatives in the $+$ and $-$
microlocalizations is controlled by the operators $\bar LL$ and $L\bar  
L$, respectively. In the $0$ microlocalization the
gain of derivatives is controlled by both $\bar LL$  and $L\bar L$  
independently. Propositions 4 and 5 give {\it a priori\/} estimates
for $E_k$ in the $0$ and $-$ microlocalizations, respectively.  
Proposition 6 gives these estimates for the $+$
microlocalization. Here we show how to go from the {\it a priori\/} 
estimates to  
hypoellipticity. In particular we prove that $E_0$ is
hypoelliptic and that $E_k$ is hypoelliptic on open sets on which  
$z\neq0$ and that the $0$ and $-$ microlocalizations of the
operators $E_k$ are hypoelliptic.

\begin{prop} If $u$ is a distribution such that for some open set  
$V\subset\mathbb R^3$ the restriction of $E_0u$ to $V$ is
in $C^\infty(V)$ then the restriction of $u$ to $V$ is also in  
$C^\infty(U)$. More precisely{\rm ,} if $E_0u\in H_{\rm loc}^s(V)$ then
$u\in H_{\rm loc}^{s+1}(V)$.
\end{prop}

\Proof  Assuming that $E_0u\in H_{\rm loc}^s(V)$, it suffices to show  
that any $P\in V$ has a neighborhood $U\subset V$ such
that for any $\varrho\in C_0^\infty(U)$ we have $\varrho u\in  
H^{s+1}(\mathbb R^3)$. Without loss of generality we may
assume that $P=0$. Now choose neighborhoods $U$ and $U'$ of $P$ such  
that $\bar U\subset U'$ and $|z|\le a$ on $U'$, as in
Proposition 4. Let $\varrho\in C_0^\infty(U)$, let $\varrho'\in  
C_0^\infty(U')$ with $\varrho'=1$ in a neighborhood of the
support of $\varrho$, and let $\theta\in C_0^\infty(\mathbb R^3)$ such  
that $\theta=1$ on a neighborhood of~$\bar U'$. Since
$u$ is a distribution there exists an $s_0\in\mathbb R$ such that  
$\theta u\in H^{-s_0}(\mathbb R^3)$. Then, choosing $\gamma^+$,
$\gamma^0$, and $\gamma^-$ such that $\gamma^++\gamma^0+\gamma^-\ge  
{\rm const.}>0$ and combining Propositions 4, 5, and 6 we obtain
the {\it a priori\/} estimate
$$
\|\varrho u\|^2_{s+1}+\|\varrho Lu\|^2_{s+\frac{1}{2}}+\|\varrho \bar  
Lu\|^2_{s+\frac{1}{2}}\le
C(\|\varrho'E_0u\|^2_s+\|u\|^2_{-s_0}),
$$
for all $u\in C^\infty(\mathbb R^3)$. Let $\chi\in C_0^\infty(\mathbb  
R^3)$ with $\chi(0)=1$. For $\delta>0$ we define the
smoothing operator $S_\delta$ by $\mathcal F(S_\delta  
u)(\xi)=\chi(\delta\xi)\hat u(\xi)$. The important facts are that:
\begin{enumerate}
\item[1.] If $\delta>0$ then for any distribution $v$ the function  
$S_\delta v\in C^\infty(\mathbb R)$.
\item[2.] If $v$ is a distribution and if $\|S_\delta v\|_s$ is bounded  
independently of $\delta$ then $v\in H^s(\mathbb R^3)$.
\item[3.] If $v\in H^s(\mathbb R^3)$ then $\lim_{\delta\to0}\|S_\delta  
v-v\|_s=0$.
\item[4.] For $\delta\ge0$ the operator $S_\delta$ is a pseudodifferential  
operator which is uniformly of order zero.
\end{enumerate}
Replacing $u$ by $S_\delta\theta u$ in Lemma 6 and in the proofs of  
Propositions 4, 5, and 6 and using item 4 above we obtain
$$
\|S_\delta\varrho u\|^2_{s+1}\le C(\|S_\delta\varrho'E_0u\|^2_s+
\|\tilde S_\delta\varrho'u\|^2_{s+\frac{1}{2}}+\|\tilde  
S_\delta\tilde\theta u\|^2_{-s_0}),
$$
where $\tilde S_\delta$ has the symbol $\tilde\chi(\delta\xi)$ with  
$\tilde\chi=1$ in a neighborhood of the support of
$\chi$. Choose $m$ so that $-s_0\ge s+1-m$, then substituting $s+1-m+j$  
for $s$ above we obtain, by induction on $j$, that
$\|S_\delta\varrho u\|^2_{s+1}$ is bounded independently of $\delta$.  
Hence $\varrho u\in H^{s+1}(\mathbb R^3)$ thus
concluding the proof.

  Next we will show that in any region in which $z\neq0$ the  
operator $E_k$ is hypoelliptic with a gain of one
derivative.
\begin{prop}If $V\subset\mathbb R^3$ is an open set{\rm ,} with the property  
that $z\neq0$ on
$V${\rm ,} and if $u$ is a distribution such that the restriction of $E_ku$  
to $V$ is in $C^\infty(V)${\rm ,} then the
restriction of $u$ to $V$ is also in $C^\infty(U)$. More precisely{\rm ,} if  
$E_ku\in H_{\rm loc}^s(V)$ then $u\in H_{\rm loc}^{s+1}(V)$.
\end{prop}

\Proof  Let $P\in V$ then $P=(\alpha,t_0)$ with $\alpha\neq0$. Let  
$U$ be a neighborhood of $P$ such that on $U$ we
have $|z-\alpha|<a$, where $a$ is chosen as in Lemma 3, and also such  
that on $U$ we have $|z|\ge b>0$. Then
$$
\|Lu\|^2\le b^{-2k}\|\bar z^kLu\|^2,
$$
for all $u\in C^\infty_0(U)$. Hence Propositions 4, 5, and 6 hold with  
$\gamma$ replaced by $\gamma_\alpha$. The proof is then
concluded using the same argument as above, replacing $S_\delta$ with  
$S_{\alpha,\delta}$, which is  defined by $\mathcal
F_\alpha(S_{\alpha,\delta}u)(\xi)=\chi(\delta\xi)\mathcal F_\alpha  
u(\xi)$.

 Now we prove microlocal hypoellipticity in the $0$ and $-$  
microlocalizations.

\begin{prop} Given neighborhoods of the origin $U$ and $U'$ with $\bar  
U\subset U'$ and $|z|\le a$ on $U'${\rm ,} where $a$ is
sufficiently small as in Lemma {\rm 3,} suppose that $\varrho\in  
C_0^\infty(U)$ and $\varrho'\in C_0^\infty(U')$ with
$\varrho'=1$ on a neighborhood of $\bar U$. Further suppose that  
$\gamma^0\in\mathcal G^0$. Then{\rm ,} given $s \in\mathbb R${\rm ,} if $u$
is a distribution such that
$\varrho'E_ku\in H^s(\mathbb R^3)$ then $\varrho\Gamma^0 u\in  
H^{s+2}(\mathbb R^3)$.
\end{prop}

\Proof  The proof consists of proving the following estimate
$$
\|S_\delta\varrho\Gamma^0u\|^2_{s+2}\le  
C(\varrho'E_ku\|^2_s+\|u\|^2_{-s_0}). 
$$
Its proof is exactly analogous to the proof of Proposition 4. Replacing  
$u$ by $S_\delta u$ the same proof as of Lemma 6 using
$XS_\delta$ instead of $X$ gives
\begin{align*}
\|S_\delta\varrho\Gamma^0Xu\|^2_{s+1}&=(\Lambda^sS_\delta\varrho\Gamma^0 
X^*Xu,\Lambda^{s+2}S_\delta\varrho\Gamma^0u)\\[4pt]
&\quad+
O(\|\tilde S_\delta\varrho_1\Gamma^0_1u\|^2_{s+1}+
\|S_\delta\varrho\Gamma^0u\|^2_{s+2}\|\tilde  
S_\delta\varrho_1\Gamma^0_1u\|^2_{s+1}\\[4pt]
&\quad+\|\tilde S_\delta\varrho_1\Gamma^0_1Xu\|^2_s+
\|\tilde  
S_\delta\varrho_1\Gamma^0_1u\|_{s+1}\|\varrho\Gamma^0Xu\|_{s+1}+\|u\|^2_ 
{-\infty}).
\end{align*}
The argument then proceeds exactly as in Proposition 4 and shows that\break
$\|S_\delta\varrho\Gamma^0u\|^2_{s+2}$ is bounded
independently of $\delta$ completing the proof.

  For the $-$ microlocalization we the following result follows  
from an argument entirely analogous to \pagebreak the above
proposition.
\begin{prop} \hskip-4pt Given neighborhoods of the origin $U$ and $U'$ with $\bar  
U\!\subset\! U'$ and $|z|\le a$ on $U'${\rm ,} where $a$ is
sufficiently small as in Lemma {\rm 3}. Suppose that $\varrho\in  
C_0^\infty(U)$ and $\varrho'\in C_0^\infty(U')$ with
$\varrho'=1$ on a neighborhood of $\bar U$. Further suppose that  
$\gamma^-\in\mathcal G^0$. Then{\rm ,} given $s \in\mathbb R${\rm ,} if $u$
is a distribution such that $\varrho'E_ku\in H^s(\mathbb R^3)$ then  
$\varrho\Gamma^- u\in H^{s+1}(\mathbb R^3)$.
\end{prop}

\section{The operator $E_1$: no loss, no gain}

As was shown in Section 5 the operator $E_1$ does not gain any  
derivatives. Here we will give a proof of an {\it a priori\/} estimate
which shows that it does not lose any derivatives. More precisely, the  
estimate will show that $E_1$ does not lose any
derivatives after it is proved that $E_1$ is hypoelliptic.  This will be  
done using the same estimate with an appropriate
smoothing operator in Section 14. As we have seen all the operators  
$E_k$ gain a derivative in regions where $z\neq0$ and in the
$0$ and $-$ microlocalizations. Thus the remaining case is the $+$  
microlocalization when $z=0$. Since the operators $E_k$
are invariant under translation in the $t$ direction it will suffice to  
consider neighborhoods of the origin. In this section
we will present a direct proof of the {\it a priori\/} 
estimates for $E_1$ which  
will rely on the following lemma. This proof however
cannot be adopted to prove the corresponding {\it a priori\/} estimate for the  
operator $F_1=E_1+c$ unless $c\ge0$. In fact the same
estimates will be proved when we treat the general case of $E_k$ with  
$k\ge1$. However that treatment is much more complicated
so it might be worthwhile to note this simpler proof.

  In the previous section we showed that the elliptic  
microlocalization $\Gamma^0u$ is smooth whenever $E_ku$ is
smooth. Thus we do not have to keep track of just which microlocalizing  
operator in $\mathfrak G^0$ is used; in
order to simplify the calculations we will write $u^0$ instead of  
$\Gamma^0u$. Similarly, since all the commutators with
$\Gamma^+$ that arise are dominated as follows $\|[\Gamma^+,R^s]u\|\le  
C(\|\Gamma^0u\|_{s-1}+\|u\|_{-\infty}$,
we will write $u^+$ instead of $\Gamma^+$.
\begin{lemma} Given a bounded open set $U\subset\mathbb R^3$ there  
exists $C>0$ such that
$$
\|u\|^2\le C(\|\bar zLu\|^2+\|\bar Lu\|^2),
$$
for all $u\in C_0^\infty(U)$.
\end{lemma}

\Proof   If $u\in C_0^\infty(U)$ we have
\begin{align*}
\|u\|^2&=(L(z)u,u)=-(zLu,u)-(zu,\bar Lu)\le \|\bar  
zLu\|\|u\|+\|zu\|\|\bar Lu\|\\
&\le {\rm s.c.}\|u\|^2+{\rm l.c.}(\|\bar zLu\|^2+\|\bar Lu\|^2).
\end{align*}
Absorbing the first term on the right into the left-hand side completes  
the proof.
 
The other estimate we will use here is given in Lemma 5 with  
$\alpha=0$, namely
\begin{equation} \label{plus}
\|\bar Lu^+\|^2+\|u^+\|_{\frac{1}{2}}^2\le
C(\|Lu^+\|^2+\|u^+\|^2+\|u\|_{-\infty}^2),
\end{equation}
for all $u\in C_0^\infty(U)$.
\begin{prop}Let $U$ be a bounded neighborhood of the origin such that  
$|z|\le a$ on $U${\rm ,} let
$\varrho,\varrho'\in C_0^\infty(U)$ with $\varrho'=1$ in a neighborhood  
of the support of $\varrho$. Then{\rm ,} given
$s,s_0\in\mathbb R$ there exists $C=C(\varrho,\varrho',s,s_0)$ such that
$$
\|\Psi^{s+\frac{1}{2}}\varrho u^+\|\le  
C(\|\Psi^{s+\frac{1}{2}}\varrho'E_1u\|+\|\Psi^s\varrho'u\|+\|u\|_{- 
s_o}),
$$
for all $u\in C_0^\infty(\mathbb R^3)$.
\end{prop}

\Proof  We assume that $u\in C_0^\infty(\mathbb R^3)$ and replace  
$u$ in (\ref{plus}) by $\varrho'\bar z\Psi^su$.  Then,
following the method of Proposition 4, we get
\begin{align*}
\|\varrho'\bar z\Psi^su^+\|_{\frac{1}{2}}^2&\le
C(\|\bar zL\varrho'\Psi^su^+\|^2+\|\bar L\varrho'\Psi^s  
u^+\|^2+\|\varrho''\Psi^su^+\|^2+\|u\|_{-\infty}^2)\\
&\le  
C(|(\varrho'\Psi^s(E_ku)^+,\varrho\Psi^su^+)|+\|\varrho''u^0\|_s^2+\| 
\varrho''\Psi^su^+\|^2+\|u\|_{-\infty}^2)\\
&\le C(\|\varrho'E_ku\|^2_s+\|\varrho''u\|_s^2+\|u\|_{-\infty}^2).
\end{align*}
Next, we replace $u$ by $\varrho\Psi^{s+\frac{1}{2}}u^+$ in Lemma 8  
and, with the use of Lemma 1 and the fact that
$$
\|\varrho'\bar  
z\Psi^su^+\|_{\frac{1}{2}}^2=\|z\varrho'\Psi^{s+\frac{1}{2}}u^+\|^2+O(\| 
u^0\|^2_{s-\frac{1}{2}}
+\|u\|_{-\infty}^2),
$$
we obtain
\begin{align*}
\|\varrho\Psi^{s+\frac{1}{2}}u^+\|^2&\le C(\|\bar  
zL\varrho\Psi^{s+\frac{1}{2}}u^+\|^2+
\|\bar  
L\varrho\Psi^{s+\frac{1}{2}}u^+\|^2+\|u^0\|^2_s+\|u\|^2_{-\infty})\\
&\le C(\|\varrho  
E_ku\|^2_{s+\frac{1}{2}}\|^2+\|L(\varrho)\Psi^{s+\frac{1}{2}}u^+\|^2\\
&\phantom{\le C(}+
\|\bar L(\varrho)\Psi^{s+\frac{1}{2}}u^+\|^2  
+\|\varrho'u^0\|^2_{s+\frac{1}{2}}+\|u\|_{-\infty}^2)\\
&\le C(\|\varrho  
E_ku\|^2_{s+\frac{1}{2}}\|^2+\|z\varrho'\Psi^{s+\frac{1}{2}}u^+\|^2+
\|\varrho'u^0\|^2_{s-\frac{1}{2}}+\|u\|_{-\infty}^2)\\
&\le C(\|\varrho' E_ku\|^2_{s+\frac{1}{2}}\|^2  
+\|\varrho''u\|_s^2+\|u\|_{-\infty}^2).
\end{align*}
Then, redefining $\varrho'$ and $\varrho''$, we conclude the proof.

\section{Estimates of $\varrho\bar Lu^+$ and of $\varrho L\bar Lu^+$}

\medskip In this section we begin to prove the {\it a priori\/} estimates for  
the operators $E_k$ with $k\ge1$. These will
be derived from the estimate (\ref{plus}) and the estimates in the $0$  
microlocalization. The main difficulty is the
localization in space; one cannot have a term with the cutoff function  
$\varrho$ between $u$ and $L$, or $\bar L$, unless the
term also contains suitable powers of $z$ and $\bar z$. Substituting  
$\varrho\Psi^s\bar Lu$ for $u$ in (\ref{plus}) we have
$$
\|\bar L\varrho\Psi^s\bar Lu^+\|^2+\|\varrho\Psi^s\bar  
Lu^+\|_{\frac{1}{2}}^2\le
C(\|L\varrho\Psi^s\bar Lu^+\|^2+\|\varrho\Psi^s\bar  
Lu^+\|^2+\|u\|_{-\infty}^2),
$$
so that,
\begin{align*}
\|\varrho&\Psi^sL\bar Lu^+\|^2+\|\varrho\Psi^s\bar  
L^2u^+\|^2+\|\varrho\Psi^{s+\frac{1}{2}}\bar Lu^+\|\\
&\le C(\|\varrho\Psi^sL\bar Lu^+\|^2+\|\varrho'\Psi^s\bar  
Lu^+\|^2+\|\varrho'u^0\|^2_{s+1}+\|u\|_{-\infty}^2)\\
&\le C(|(\varrho\Psi^s\bar LL\bar Lu^+,\varrho\Psi^s\bar  
Lu^+)|+\|\varrho'\Psi^s\bar Lu^+\|^2+
\|\varrho''E_ku\|^2_{s-1}+\|u\|_{-\infty}^2).
\end{align*}
Since $\bar LL\bar L=-\bar LE_k-\bar L^2|z|^{2k}L$, we have
\begin{eqnarray*}
&&\hskip-12pt |(\varrho\Psi^s\bar LL\bar Lu^+,\varrho\Psi^s\bar Lu^+)|\\
&&\le  
C(|(\varrho\Psi^s\bar LE_ku^+,\varrho\Psi^s\bar Lu^+)|+
|(\varrho\Psi^s\bar L^2|z|^{2k}Lu^+,\varrho\Psi^s\bar Lu^+)|\\
&&\le {\rm l.c.}\|\varrho'E_ku\|^2_s+{\rm s.c.}\|\varrho\Psi^sL\bar Lu^+\|^2+
C|(\varrho\Psi^s\bar L^2|z|^{2k}Lu^+,\varrho\Psi^s\bar Lu^+)|.
\end{eqnarray*}
Then, to estimate $|(\varrho\Psi^s\bar  
L^2|z|^{2k}Lu^+,\varrho\Psi^s\bar Lu^+)|$, we have
\begin{align*}
\bar L^2|z|^{2k}L&=-k\bar Lz^k\bar z^{k-1}L+\bar L|z|^{2k}\bar LL\\
&=-k^2\bar L|z|^{2(k-1)}+\bar LLz^k\bar z^{k-1}-2kz^k\bar z^{k-1}T+\bar  
L|z|^{2k}L\bar L-2|z|^{2k}T\bar L\\
&=-k^2\bar L|z|^{2(k-1)}-4kz^k\bar z^{k-1}T+k(k-1)Lz^k\bar  
z^{k-2}+kLz^k\bar z^{k-1}\bar L\\*
&\quad+\bar L|z|^{2k}L\bar L-2|z|^{2k}T\bar L,
\end{align*}
and, using integration by parts, we get
\begin{align*}
|(\varrho\Psi^s\bar L|z|^{2(k-1)}u^+,\varrho\Psi^s\bar Lu^+)|&\le  
{\rm l.c.}\|z^{2(k-1)}\varrho\Psi^su^+\|^2+\mathcal E_1,\\[4pt]
|(\varrho\Psi^sz^k\bar z^{k-1}Tu^+,\varrho\Psi^s\bar Lu^+)|&\le  
{\rm l.c.}\|z^{2k-1}\varrho\Psi^{s+\frac{1}{2}}u^+\|^2+\mathcal E_2,\\[4pt]
(k-1)|(\varrho\Psi^sLz^k\bar z^{k-2}u^+,\varrho\Psi^s\bar Lu^+)|&\le  
(k-1)({\rm l.c.}\|z^{2(k-1)}\varrho\Psi^su^+\|^2+\mathcal E_3),\\[4pt]
 |(\varrho\Psi^sLz^k\bar z^{k-1}\bar Lu^+,\varrho\Psi^s\bar Lu^+)|&\le  
{\rm l.c.}\|z^{2k-1}\varrho\Psi^su^+\|^2+\mathcal E_4,\\[4pt]
 |(\varrho\Psi^s\bar L|z|^{2k}L\bar Lu^+,\varrho\Psi^s\bar Lu^+)|&\le  
\mathcal E_4,\\[4pt]
\noalign{\noindent \rm and} 
|(\varrho\Psi^s|z|^{2k}T\bar Lu^+,\varrho\Psi^s\bar Lu^+)|&\le \mathcal  
E_2,
\end{align*}
where
\begin{align*}
\mathcal E_1&\sim\|\varrho'u^0\|^2_s+\|\varrho'\Psi^s\bar  
Lu^+\|^2+\|u\|^2_{-\infty},\\[4pt]
 \mathcal E_2&\sim {\rm s.c.}\|\varrho\Psi^{s+\frac{1}{2}}\bar  
Lu^+\|^2+\mathcal E_1,\\[4pt]
 \mathcal E_3&\sim {\rm s.c.}\|\varrho\Psi^s\bar L^2u^+\|^2+\mathcal E_1,\\[4pt]
\noalign{\noindent\rm and} 
 \mathcal E_4&\sim {\rm s.c.}\|\varrho\Psi^sL\bar Lu^+\|^2+\mathcal E_1. 
\end{align*} 

\vglue-23pt
\begin{align*}
|(\Psi^s\varrho\bar L|z|^{2(k-1)}u^+,\Psi^s\varrho \bar Lu^+)|&\le  
C(\|z^{2k-2}\Psi^s\varrho u^+\|^2
+\mathcal E_2),\\[4pt]
|(\Psi^s\varrho\bar Lz^{k-1}\bar z^k\bar Lu^+,\Psi^s\varrho \bar  
Lu^+)|&\le C\mathcal E_1\\[4pt]
|(\Psi^sz^k\bar z^{k-1}\varrho Tu^+,\Psi^s\varrho \bar Lu^+)|&\le  
C(\|z^{2k-1}\Psi^{s+\frac{1}{2}}\varrho u^+\|^2+
\mathcal E_1+\mathcal E_3)\\[4pt]
|(\Psi^s|z|^{2k}T\varrho\bar Lu^+,\Psi^s\varrho \bar Lu^+)|&\le  
C(\mathcal E_1+\mathcal E_2)\\[4pt]
|(\Psi^sLz^k\bar z^{k-2}\varrho u^+,\Psi^s\varrho \bar Lu^+)|&\le  
C(\|z^{2k-2}\Psi^s\varrho u^+\|^2+
\mathcal E_1+\mathcal E_4)\\[4pt]
|(\Psi^sLz^k\bar z^{k-1}\varrho\bar Lu^+,\Psi^s\varrho \bar Lu^+)|&\le  
C(\mathcal E_1+\mathcal E_4),\\[4pt]
\noalign{\noindent \rm and} 
 |(\Psi^sL|z|^{2k}\varrho\bar L^2u^+,\Psi^s\varrho \bar Lu^+)|&\le  
C(\mathcal E_1+\mathcal E_4).
\end{align*}
Again, let $\{\varrho_i\}$ be a sequence of cutoff functions as defined  
in Section 2. Then substituting $\varrho_{i}$ for
$\varrho$, $s-\frac{i-1}{2}$ for $s$, and $\varrho_{i+1}$ for  
$\varrho'$, we get
\begin{align*}
&\hskip-12pt\|\varrho_i\Psi^{s-\frac{i-1}{2}}L\bar  
Lu^+\|^2+\|\varrho_i\Psi^{s+1-\frac{i}{2}}\bar Lu^+\|^2\\[4pt]
&\quad\le
C(\|\varrho'E_ku\|^2_{s-\frac{i-1}{2}}+\|z^{2k-2}\Psi^{s-\frac{i- 
1}{2}}\varrho_i u^+\|^2\\[4pt]
&\qquad+\|z^{2k-1}\Psi^{s+1-\frac{i}{2}}\varrho_i u^+\|^2  
+\|\varrho_{i+1}\Psi^{s-\frac{i-1}{2}}\bar Lu^+\|^2
+\|u\|^2_{-\infty}).
\end{align*}
Then we obtain the following, by substituting these inequalities into  
each other for successive $i$
\begin{eqnarray*}
&&\|\varrho\Psi^{s}L\bar Lu^+\|^2+\|\varrho\Psi^{s+\frac{1}{2}}\bar  
Lu^+\|^2\\[4pt]
&&  \le
C\Big(\sum_{i=1}^N\big(\|\varrho_i\Psi^{s-\frac{i-1}{2}}\bar  
LE_ku^+\|^2\!+\|z^{2k-2}\Psi^{s-\frac{i-1}{2}}\varrho_i
u^+\|^2+\|z^{2k-1}\Psi^{s+1-\frac{i}{2}}\varrho_i u^+\|^2\big)  
\\[4pt]
&&\qquad+\|\varrho_{N+1}\Psi^{s-\frac{N-1}{2}}\bar Lu^+\|^2
+\|\varrho_NE_ku\|_{s-1}^2+\|u\|^2_{-\infty}\Big).
\end{eqnarray*}
Given $s_o$ we choose $N>2(s-s_o)+1$ then we obtain the following  
estimate which will be repeatedly used in establishing
the {\it a priori\/} estimates for $E_k$
\begin{align}\label{central}
\|\varrho\Psi^sL\bar Lu^+\|^2&+\|\varrho\Psi^{s+\frac{1}{2}}\bar  
Lu^+\|^2\le
C(\|\varrho'E_ku^+\|^2_s+\|z^{2k-2}\Psi^s\varrho' u^+\|^2\\[4pt]
&\notag+\|z^{2k-1}\Psi^{s+\frac{1}{2}}\varrho' u^+\|^2 +\|u\|^2_{-s_o}).
\end{align}

\section{Estimates of $\|z^j\Psi^{s+ja}\varrho u^+\|$}

  \begin{lemma} If $a>0$ then for $m\in\mathbb Z^+$ and a small  
constant\/ ${\rm s.c.}$ there exists a constant\/ ${\rm l.c.}$ such that
\begin{eqnarray*}
\sum_{j=1}^{m-1}\|z^j\Psi^{s+ja}\varrho u^+\|^2&\le & 
{\rm l.c.}\|z^m\Psi^{s+ma}\varrho u^+\|^2+{\rm s.c.}\|\Psi^s\varrho u\|^2
\\*
&&+C(\|\varrho'u^0\|^2_{s+(m-1)a-1}+\|u\|^2_{-\infty}),\\*
\noalign{\noindent for  all $u\in C^\infty(U)$.}
\end{eqnarray*}
\end{lemma}

\vskip-12pt
{\it Proof}.  For $m=2$ we have
\begin{align*}
\|z\Psi^{s+a}\varrho u^+\|^2&=(|z|^2\Psi^{s+2a}\varrho  
u^+,\Psi^{s}\varrho u^+)+O(\|\varrho'u^0\|^2_{s+a-1}
+\|u\|^2_{-\infty})\\
&\le {\rm l.c.}\|z^2\Psi^{s+2a}\varrho u^+\|^2+{\rm s.c.}\|\Psi^{s}\varrho  
u^+\|^2+C(\|\varrho'u^0\|^2_{s+a-1}+\|u\|^2_{-\infty}).
\end{align*}
For $m>2$ we assume
\begin{eqnarray*}
\sum_{j=1}^{m-2}\|z^j\Psi^{s+ja}\varrho u^+\|^2&\le&
{\rm l.c.}\|z^{m-1}\Psi^{s+(m-1)a}\varrho u^+\|^2+{\rm s.c.}\|\Psi^s\varrho  
u\|^2\\
&&+C(\|\varrho u^0\|^2_{s+(m-2)a-1}+\|u\|^2_{-\infty}),
\end{eqnarray*}
and we have
\begin{eqnarray*}
\|z^{m-1}\Psi^{s+(m-1)a}\varrho u^+\|^2&=&(z^{m}\bar z\Psi^{s+ma}\varrho  
u^+,z^{m-2}\Psi^{s+(m-2)a}\varrho u^+)\\
&&+O(\|\varrho
u^0\|^2_{s+(m-1)a-1}+\|u\|^2_{-\infty})\\
& \le& {\rm l.c.}\|z^{m}\Psi^{s+ma}\varrho  
u^+\|^2+{\rm s.c.}\|z^{m-2}\Psi^{s+(m-2)a}\varrho u^+\|^2\\
&&+
C(\|\varrho u^0\|^2_{s+(m-1)a-1}+\|u\|^2_{-\infty}).
\end{eqnarray*}
Adding this to the above and absorbing the term multiplied by ${\rm s.c.}$ in  
the right-hand side we conclude the proof.

 \begin{lemma} If $0<j<m$ and if $\frac{mA}{j}<B$ then for any  
${\rm s.c.}$ and any $N$ there exists $C_N$ such that
\begin{eqnarray*}
\|z^j\Psi^{s+A}\varrho u^+\|^2&\le& {\rm s.c.}(\|z^m\Psi^{s+B}\varrho  
u^+\|^2+\|\Psi^{s}\varrho u^+\|^2)\\
&&+C(\|\varrho u^0\|^2_{s+B-1}+C_N\|u^+\|^2_{-N}+C(\|\varrho  
u^0\|^2_{s+B-1}+\|u\|^2_{-\infty}),
\end{eqnarray*}
for all $u\in C^\infty_0(U)$.
\end{lemma}

\Proof  With $a=\frac{A}{j}$ we have
\begin{eqnarray*}
\|z^j\Psi^{s+A}\varrho u^+\|^2&\le& {\rm l.c.}\|z^m\Psi^{s+ma}\varrho  
u^+\|^2+s.c\|\Psi^{s}\varrho u^+\|^2\\
&&+
C(\|\varrho u^0\|^2_{s+ma-1}+\|u\|^2_{-\infty}).
\end{eqnarray*}
Since $ma=\frac{mA}{j}<B$,
$$
\psi^{s+ma}(\xi)\le {\rm s.c.}\psi^{s+B}(\xi)+{\rm l.c.}(1+|\xi|^2)^{-\frac{N}{2}}.
$$
Then
\begin{eqnarray*}
\|z^m \Psi^{s+ma}\varrho u^+\|^2&=&\|\Psi^{s+ma}z^m\varrho  
u^+\|^2+O(\|u^0\|^2_{s+ma-1}+\|u\|^2_{-\infty})\\
&\le& {\rm s.c.}\|\Psi^{s+B}z^m\varrho u^+\|^2+C_N\|u^+\|^2_{-N}\\
&&+O(\|\varrho  
u^0\|^2_{s+ma-1}+\|u\|^2_{-\infty})\\
&\le &{\rm s.c.}\|z^m\Psi^{s+B}\varrho u^+\|^2+C_N\|u^+\|^2_{-N}\\
&&+O(\|\varrho  
u^0\|^2_{s+B-1}+\|u\|^2_{-\infty}).
\end{eqnarray*}
Combining with the above we conclude  the proof of the lemma.

\begin{lemma} If $\sigma=\frac{1}{2k}$ and if $1\le j\le k$ then
$$
\|z^j\Psi^{s+j\sigma}\varrho u^+\|^2\le  
C(\|\varrho'E_ku\|^2_s+\|\varrho'u\|^2_s+\|u\|^2_{-\infty}),
$$
for all $u\in C^\infty_0(U)$.
\end{lemma}

\Proof  First note that
$$
\|\varrho z^k\Psi^su^+\|^2\le  
C(\|\varrho'E_ku\|^2_{s-\frac{1}{2}}+\|\varrho'\Psi^{s- 
\frac{1}{2}}u^+\|^2+\|u\|^2_{-\infty}).
$$
Then, replacing $s$ by $s+k\sigma$, since $k\sigma-\frac{1}{2}=0$, we  
have
\begin{align*}
\|z^j\Psi^{s+j\sigma}\varrho u^+\|^2&\le C(\|z^k\Psi^{s+k\sigma}\varrho  
u^+\|^2+\|\Psi^s\varrho
u\|^2\\
&\quad+C\|\varrho'u^0\|^2_{s+(k-1)a-1} +\|u\|^2_{-\infty})\\
&\le C(\|\bar z^kL\Psi^s\varrho u^+\|^2+\|\bar L\Psi^s\varrho  
u^+\|^2+\|\varrho'u^+\|^2_s+\|u\|^2_{-\infty})\\
&\le C(\|\varrho'E_ku\|^2_s+\|\varrho'u\|^2_s+\|u\|^2_{-\infty}).
\end{align*}

\section{Estimate of $\|\varrho\Psi^{s+\sigma}u^+\|$}

\begin{lemma} There exists a $C>0$ such that
$$
\|\varrho\Psi^{s+\sigma}u^+\|\le  
C(\|\varrho'E_ku\|^2_{s+\sigma+k-1}+\|\varrho'u\|^2_s+\|u\|^2_{- 
\infty}),
$$
for all $u\in C^\infty_0(U)$.
\end{lemma}

\Proof 
\begin{align*}
\|\Psi^{s+\sigma}\varrho u^+\|^2&=(L(z)\Psi^{s+\sigma}\varrho  
u^+,\Psi^{s+\sigma}\varrho u^+)\\
&=-(zL\Psi^{s+\sigma}\varrho u^+,\Psi^{s+\sigma}\varrho  
u^+)-(z\Psi^{s+\sigma}\varrho u^+,\bar L\Psi^{s+\sigma}\varrho u^+)\\
&\le {\rm l.c.}\|z\Psi^{s+\sigma}\varrho  
Lu^+\|^2+C\|\Psi^{s+\sigma}\varrho\bar L u^+\|^2+``error",
\end{align*}
where,
$$
``{\rm error}"\le {\rm s.c.}\|\Psi^{s+\sigma}\varrho  
u^+\|^2+C(\|z\Psi^{s+\sigma}\varrho' u^+\|^2+\|\varrho  
u^0\|^2_{s+\sigma}+\|u\|^2_{-\infty}).
$$
In the estimate of the ``error'' the first term on the right gets  
absorbed and the other terms are estimated as follows.
$$
\|z\Psi^{s+\sigma}\varrho' u^+\|^2\le C(\|\varrho'  
Eu\|^2_s+\|\varrho'u\|^2_s+\|u\|^2_{-\infty}).
$$
The third term, which is microlocalized in the elliptic region, is  
estimated by
$$
\|\varrho u^0\|^2_{s+\sigma}\le C(\|\varrho  
Eu\|^2_{s+\sigma-2}+\|\varrho'u\|^2_s).
$$
Hence we get
\begin{eqnarray*}
\|\Psi^{s+\sigma}\varrho u^+\|^2&\le&C\Big(\|z\Psi^{s+\sigma}\varrho  
Lu^+\|^2+\|\Psi^{s+\sigma}\varrho\bar Lu^+\|^2\\
&&\phantom{C\big(}+\|\varrho Eu\|^2_s
+\|\varrho'u\|^2_s+\|u\|^2_{-\infty}\Big).
\end{eqnarray*}
 From (\ref{central}) we have
$$
\|\Psi^{s+\sigma}\varrho\bar Lu^+\|^2\le C(\|\varrho  
Eu\|^2_{s+\sigma-\frac{1}{2}}+\|\varrho'u\|^2_s+\|u\|^2_{-\infty}).
$$
So the term that remains to be estimated is $\|z\Psi^{s+\sigma}\varrho  
Lu^+\|^2$, and  we have
\begin{eqnarray*}
\|z\Psi^{s+\sigma}\varrho  
Lu^+\|^2&=&(|z|^2\Psi^{s+\sigma+\frac{1}{2}}\varrho  
Lu^+,\Psi^{s+\sigma-\frac{1}{2}}\varrho Lu^+)
\\
&&+O(\|u^0\|^2_{s+\sigma-2}+\|u\|^2_{-\infty})\\
&\le& {\rm l.c.}\|z^2\Psi^{s+\sigma+\frac{1}{2}}\varrho  
Lu^+\|^2+{\rm s.c.}\|\Psi^{s+\sigma-\frac{1}{2}}\varrho Lu^+\|^2
\\
&&+O(\|u^0\|^2_{s+\sigma-2}+\|u\|^2_{-\infty})
\end{eqnarray*}
and
\begin{eqnarray*}
\|\Psi^{s+\sigma-\frac{1}{2}}\varrho Lu^+\|^2&=&(\Psi^{s+\sigma}\varrho  
Lu^+,\Psi^{s+\sigma-1}\varrho Lu^+)
\\
&=&-(\Psi^{s+\sigma}\varrho u^+,\bar L\Psi^{s+\sigma-1}\varrho  
Lu^+)+\mathcal E_1\\
&=&-(\Psi^{s+\sigma}\varrho u^+,[\bar L,\Psi^{s+\sigma-1}\varrho  
L]u^+)\\
&&-(\Psi^{s+\sigma}\varrho u^+,\Psi^{s+\sigma-1}\varrho L\bar
Lu^+)+\mathcal E_1\\
&\le& C(\|\Psi^{s+\sigma}\varrho u^+\|^2+\|[\bar  
L,\Psi^{s+\sigma-1}\varrho L]u^+\|^2 \\
&&+\|\Psi^{s+\sigma-1}\varrho L\bar
Lu^+\|^2+\mathcal E_1.
\end{eqnarray*}
The second term is estimated as follows
\begin{eqnarray*}
[\bar L,\Psi^{s+\sigma-1}\varrho L]u^+&=&[\bar  
L,\Psi^{s+\sigma-1}]\varrho Lu^++\Psi^{s+\sigma-1}\bar L(\varrho)Lu^+
\\
&&-2\Psi^{s+\sigma-1}\varrho Tu^++\Psi^{s+\sigma-1}\varrho L\bar Lu^+
\end{eqnarray*}
so that
\begin{eqnarray*}
\|[\bar L,\Psi^{s+\sigma-1}]\varrho Lu^+\|^2&\le  &
C(\|\Psi^{s+\sigma-1}\Gamma^0\varrho u^+\|^2+\|u\|^2_{-\infty})\\
&\le  &
C(\|\varrho'u^0\|^2_{s+\sigma}
+\|u\|^2_{-\infty}),
\end{eqnarray*}
and
\begin{multline*}
\|\Psi^{s+\sigma-1}\bar L(\varrho)Lu^+\|^2+\|\Psi^{s+\sigma-1}\varrho  
Tu^+\|^2\\
\le C(\|z\Psi^{s+\sigma}\varrho'u^+\|^2
+\|\Psi^{s+\sigma}\varrho u\|^2)+\mathcal E_2.
\end{multline*}
Furthermore we have
$$
\|\Psi^{s+\sigma-1}\varrho L\bar Lu^+\|^2\le C\|\varrho'E_k    
u\|^2_{s+\sigma-\frac{1}{2}}+\mathcal E_3.
$$
The terms $\mathcal E$ are bounded as follows
$$
\mathcal E_1\le  
C(\|u^0\|^2_{s+\sigma}+\|z\Psi^{s+\sigma}\varrho'u^+\|^2+\|u\|^2_{- 
\infty}).
$$
By Lemma 10 we get
\begin{eqnarray*}
\mathcal E_1&\le&  
C(\|\varrho'E_ku\|^2_s+\|\varrho'u\|^2_s+\|u\|^2_{-\infty}),
\\[4pt]
\mathcal E_2&\le& C(\varrho'u\|^2_{s+\sigma-1}+\mathcal E_1)\le  
C'\mathcal E_1,
\end{eqnarray*}
and
$$
\mathcal E_3\le C(\|z^{2k-1}\varrho u^+\|^2_s+\|z^{2k-2}\varrho  
u^+\|^2_{s-\frac{1}{2}}+\mathcal E_2)\le C'\mathcal E_2.
$$
Hence we have
\begin{multline*}
\|\Psi^{s+\sigma}\varrho u^+\|^2+\|z\Psi^{s+\sigma}\varrho Lu^+\|^2\\
\le  
C(\|z^2\Psi^{s+\sigma+\frac{1}{2}}\varrho Lu^+\|^2+
\|\varrho'E_ku\|^2_s+\|\varrho'u\|_s^2+\|u\|^2_{-\infty}).
\end{multline*}
To estimate the first term on the right we will use Lemma 8 as follows.
$$
\|z^2\Psi^{s+\sigma+\frac{1}{2}}\varrho Lu^+\|^2\le  
C(\|z\Psi^{s+\sigma+\frac{1}{2}}\varrho zLu^+\|^2+\|\varrho  
u^0\|^2_{s+\sigma}
+\|u\|^2_{-\infty}).
$$
We apply Lemma 8 with $a=\frac{1}{2}$, $m=k-1$, $s$ replaced by  
$s+\sigma$, and $u$ replaced by $zLu$ to obtain
\begin{align*}
\|z^2\Psi^{s+\sigma+\frac{1}{2}}\varrho Lu^+\|^2&\le  
{\rm l.c.}\|z^{k-1}\Psi^{s+\sigma+\frac{k-1}{2}}\varrho zLu^+\|^2+
{\rm s.c.}\|z\Psi^{s+\sigma}\varrho Lu^+\|^2\\[4pt]
&\quad+\|\varrho u^0\|^2_{s+\sigma+\frac{k-1}{2}}+\|u\|^2_{-\infty})\\[4pt]
&\le {\rm l.c.}\|z^k\Psi^{s+\sigma+\frac{k-1}{2}}\varrho Lu^+\|^2+  
{\rm s.c.}\|z\Psi^{s+\sigma}\varrho Lu^+\|^2\\[4pt]
&\quad+C(\|\varrho u^0\|^2_{s+\sigma+\frac{k-1}{2}}+\|u\|^2_{-\infty}).
\end{align*}
Therefore we have
\begin{align*}
&\hskip-16pt\|\Psi^{s+\sigma}\varrho u^+\|^2\\[4pt]
& \le  
C(\|z^k\Psi^{s+\sigma+\frac{k-1}{2}}\varrho  
Lu^+\|^2+\|\varrho'E_ku\|^2_{s+\sigma+\frac{k-1}{2}-2}+
\|\varrho'u\|_s^2+\|u\|^2_{-\infty})\\[4pt]
& \le C(\|\Psi^{s+\sigma+\frac{k-1}{2}}\varrho\bar  
z^kLu^+\|^2+\|\varrho'E_ku\|^2_{s+\sigma+\frac{k-1}{2}-2}+
\|\varrho'u\|_s^2+\|u\|^2_{-\infty}).
\end{align*}
Next, from Lemma 8 with $m=k$, $a=\frac{1}{2}$ and $s$ replaced by  
$s+\sigma$, we have
\begin{eqnarray*}
\|z\Psi^{s+\sigma+\frac{1}{2}}u\|^2&\le  &
{\rm l.c.}\|z^k\Psi^{s+\sigma+\frac{k}{2}}\varrho  
u^+\|^2\\[4pt]
&&+{\rm s.c.}\|\Psi^{s+\sigma}\varrho u^+\|^2+
C(\|\varrho'u^0\|^2_{s+\sigma+\frac{k-1}{2}-1}+\|u\|^2_{-\infty})\\[4pt]
&\le &C\|\Psi^{s+\sigma+\frac{k-1}{2}}\varrho\bar z^kLu^+\|^2+\mathcal  
E_1
\end{eqnarray*}
\pagebreak

\noindent
and
\begin{align*}
\|z^k\Psi^{s+\sigma+\frac{k}{2}}\varrho u^+\|^2&\le  
C\|\Psi^{s+\sigma+\frac{k-1}{2}}\varrho\bar z^kLu^+\|^2+\mathcal E_1\\
&=-C(\Psi^{s+\sigma+\frac{k-1}{2}}\varrho\bar  
L|z|^{2k}Lu^+,\Psi^{s+\sigma+\frac{k-1}{2}}\varrho u^+)\\
&\quad-2C(\Psi^{s+\sigma+\frac{k-1}{2}}\varrho\bar  
z^{k}Lu^+,\Psi^{s+\sigma+\frac{k-1}{2}}\bar z^{k+1}\mu u^+)+\mathcal E_2
\end{align*}
since
\begin{align*} \hskip-8pt
|(\Psi^{s+\sigma+\frac{k-1}{2}}\varrho\bar  
z^{k}Lu^+,\Psi^{s+\sigma+\frac{k-1}{2}}\bar z^{k+1}\mu u^+)|&\le
{\rm s.c.}\|z^{k}\Psi^{s+\sigma+\frac{k}{2}}\varrho Lu^+\|^2\\
&\quad+{\rm l.c.}\|z^{k+1}\Psi^{s+\sigma+\frac{k- 
2}{2}}\varrho'u^+\|^2+\mathcal E_2.
\end{align*}
Hence we obtain
\begin{align*}
\|z^k\Psi^{s+\sigma+\frac{k}{2}}\varrho u^+\|^2 &\le  
C\|\Psi^{s+\sigma+\frac{k-1}{2}}\varrho\bar z^kLu^+\|^2+\mathcal E_2\\
&\le
C|(\Psi^{s+\sigma+\frac{k-1}{2}}\varrho\bar  
L|z|^{2k}Lu^+,\Psi^{s+\sigma+\frac{k-1}{2}}\varrho u^+)|+\mathcal E_3\\
&\le C|(\Psi^{s+\sigma+\frac{k-1}{2}}\varrho  
E_ku^+,\Psi^{s+\sigma+\frac{k-1}{2}}\varrho u^+)|\\
&\quad+
|(\Psi^{s+\sigma+\frac{k-1}{2}}\varrho L\bar  
Lu^+,\Psi^{s+\sigma+\frac{k-1}{2}}\varrho u^+)|+\mathcal E_3\\
&\le C(\|\Psi^{s+\sigma+k-1}\varrho E_ku\|^2\\
&\quad+
|(\Psi^{s+\sigma+\frac{k-1}{2}}\varrho \bar  
Lu^+,\Psi^{s+\sigma+\frac{k-1}{2}}\bar z\mu u^+)|)+\mathcal E_4\\
&\le C(\|\Psi^{s+\sigma+k-1}\varrho  
E_ku\|^2+\|\Psi^{s+\sigma+k-1}\varrho \bar Lu^+\|^2)+\mathcal E_5\\
&\le C(\|\varrho  
E_ku\|_{s+\sigma+k-1}^2+\|z^{2k-1}\Psi^{s+\sigma+k-1}\varrho' u^+\|^2\\
&\quad +\|z^{2k-2}\Psi^{s+\sigma+k-\frac{3}{2}}\varrho'  
u^+\|^2)+\mathcal E_6.
\end{align*}
Thus, applying Lemma 9 with $m=2k-1$, $j=2k-2$, $A=k-\frac{3}{2}$,  
$B=k-1$, and $s$ replaced by $s+\sigma$, we have
$$
\frac{mA}{j}=\frac{2k-1}{2k-2}(k-\frac{3}{2})<k-1=B.
$$
Now,
$$
\|z^{2k-2}\Psi^{s+\sigma+k-\frac{3}{2}}\varrho' u^+\|^2\le  
{\rm s.c.}\|z^{2k-1}\Psi^{s+\sigma+k-1}\varrho' u^+\|^2+\mathcal E_7.
$$
Replacing $\varrho u^+$ by $\bar z^{k-1}\varrho'u^+$ and $s$ by  
$s+\frac{k-2}{2}$ we obtain
\begin{align*} &\hskip-3pt
\|z^{2k-1}\Psi^{s+\sigma+k-1}\varrho' u^+\|^2\\&\le
C(\|\Psi^{s+\sigma+k-\frac{3}{2}}\varrho'\bar z^{2k-1}Lu^+\|^2+\mathcal  
E_8\\
&\le
C(|(\Psi^{s+\sigma+k-\frac{3}{2}}\bar  
L\varrho'|z|^{2k}Lu^+,|z|^{2k-2}\Psi^{s+\sigma+k-\frac{3}{2}}\varrho'  
u)|\\
&\quad+{\rm s.c.}\|\Psi^{s+\sigma+k-1}\varrho'\bar z^{2k-1}Lu^+\|^2+  
{\rm l.c.}\|z^{2k-1}\Psi^{s+\sigma+k-2}\varrho''
u\|^2) +\mathcal E_8\\
&\le C(\|\varrho' E_ku\|^2_{s+\sigma+2k-3}+
|(\Psi^{s+\sigma+k-\frac{3}{2}}\varrho' L\bar  
Lu^+,|z|^{2k-2}\Psi^{s+\sigma+k-\frac{3}{2}}\varrho' u)|)+\mathcal E_9\\
&\le C(\|\varrho'  
E_ku\|^2_{s+\sigma+2k-3}+|(\Psi^{s+\sigma+k-2}\varrho''\bar  
Lu^+,z|z|^{2k-2}\Psi^{s+\sigma+k-1}\varrho' u)|\\
&\quad+|(\Psi^{s+\sigma+k-1}\varrho'\bar Lu^+,z^{k-1}\bar  
z^{k-2}\Psi^{s+\sigma+k-2}\varrho' u)|)+\mathcal E_9\\
&\le C(\|\varrho'  
E_ku\|^2_{s+\sigma+2k-3}+{\rm l.c.}\|\Psi^{s+\sigma+k-2}\varrho''\bar  
Lu^+\|^2+
{\rm s.c.}\|z^{2k-1}\Psi^{s+\sigma+k-1}\varrho' u\|^2\\
&\quad+{\rm s.c.}\|\Psi^{s+\sigma+k-1}\varrho'\bar  
Lu^+\|^2+{\rm l.c.}\|z^{2k-3}\Psi^{s+\sigma+k-2}\varrho' u^+\|^2)+\mathcal  
E_9.
\end{align*}
Now applying Lemma 9 as above but with $j=2k-3$ and $A=k-2$, we get
$$
\frac{mA}{j}=\frac{2k-1}{2k-3}(k-2)<k-1=B.
$$
Hence
$$
\|z^{2k-3}\Psi^{s+\sigma+k-2}\varrho' u^+\|^2\le  
{\rm s.c.}\|z^{2k-1}\Psi^{s+\sigma+k-1}\varrho' u^+\|^2+\mathcal E_9.
$$
Combining the above we obtain
$$
\|\Psi^{s+\sigma}\varrho u^+\|^2\le C\|\varrho'  
E_ku\|_{s+\sigma+k-1}^2+\mathcal E_{10}.
$$
To complete the proof of the {\it a~priori\/} estimate we will analyze the  
error terms:
\begin{align*}
&\mathcal E_1\sim\|z^k\Psi^{s+\sigma+\frac{k-1}{2}}\varrho  
u^+\|^2+\|\varrho'u^0\|^2_{s+\sigma+\frac{k-1}{2}}
+{\rm s.c.}\|\Psi^{s+\sigma}\varrho u^+\|^2+\|u\|^2_{-\infty},\\[4pt]
&\mathcal E_2\sim\mathcal  
E_1+\|z^{k+1}\Psi^{s+\sigma+\frac{k-1}{2}}\varrho'u^+\|^2,\\[4pt]
&\mathcal E_3\sim\mathcal  
E_2+\|\varrho'u^0\|^2_{s+\sigma+\frac{k}{2}}+{\rm s.c.}\|\Psi^{s+\sigma+\frac{ 
k-1}{2}}\varrho\bar z^kLu^+\|^2,\\[4pt]
&\mathcal E_4\sim\mathcal  
E_3+\|\Psi^{s+\sigma+\frac{k-1}{2}}\varrho\bar Lu^+\|^2,\\[4pt]
&\mathcal E_5\sim\mathcal E_4+\|z\Psi^{s+\sigma}\varrho'u^+\|^2\\[4pt]
&\mathcal E_6\sim\mathcal  
E_5+\|z^{2k-1}\Psi^{s+\sigma+k-\frac{3}{2}}\varrho'u^+\|^2  
+\|z^{2k-2}\Psi^{s+\sigma+k-2}\varrho'u^+\|^2\\[4pt]
&\mathcal E_7\sim\mathcal E_6+\|\varrho u^0\|^2_{s+k-2}+\|u^+\|^2_{-N}\\[4pt]
&\mathcal E_8\sim\mathcal  
E_7+\|\Psi^{s+\sigma+k-\frac{3}{2}}\varrho'\bar z^{2k}Lu^+\|^2\\[4pt]
&\mathcal E_9\sim\mathcal E_8+{\rm s.c.}\|\Psi^{s+\sigma+k-1}\varrho\bar  
z^{2k-1}Lu^+\|^2+\|z^{2k-1}\Psi^{s+\sigma+k-2}\varrho'u\|^2,\\[4pt]
\noalign{\noindent {\rm and}}
&\mathcal E_{10}\sim\mathcal E_9+{\rm s.c.}\|\Psi^{s+\sigma+k-1}\varrho\bar  
Lu^+\|^2+\|\Psi^{s+\sigma+k-2}\varrho'\bar Lu^+\|^2.
\end{align*}
The ``admissible" errors are $\|\varrho'u\|^2_s+\|u\|^2_{-\infty}$. The  
terms involving $u^0$ are all bounded by
${\rm const}.\|\varrho'E_ku\|^2_{s+\sigma+k-2}$ modulo admissible errors. The  
terms involving a small constant
${\rm s.c.}$ are absorbed in the left. The term  
$\|z\Psi^{s+\sigma}\varrho'u^+\|$ is bounded by  
${\rm const.}\|\varrho'E_ku\|^2_s$, and the remaining terms can
be bounded by a constant times $\mathcal A(s,\varrho')$, where  
$\mathcal A(s,\varrho')$ is defined by
$$
\mathcal A(s,\varrho')=\|z^k\Psi^{s+\sigma+\frac{k-1}{2}}\varrho'u^+\|^2+\|  
z^{2k-1}\Psi^{s+\sigma+k-\frac{3}{2}}\varrho'u^+\|^2
+\|\Psi^{s+\sigma+k-2}\varrho'\bar Lu^+\|^2.  
$$
Repeating the same estimates with $s$ replaced by $s-\frac{1}{2}$ we  
replace the error $\mathcal A(s,\varrho')$ by
$\mathcal A(s-\frac{1}{2},\varrho'')$. Repeating this process $2k-2$  
times (and redefining $\varrho'$) we obtain the desired
{\it a~priori\/} estimate, namely:
\begin{equation}\label{final}
\|\Psi^{s+\sigma}\varrho u^+\|^2\le C(\|\varrho'  
E_ku\|_{s+\sigma+k-1}^2+\|\varrho'u\|^2_s+\|u\|^2_{-\infty}).
\end{equation}

\section{Smoothing}

 To conclude the proof of Theorem C we will apply the above  
estimate to the smoothing of a solution. Given a
distribution solution $u$ of $E_ku=f$ with $f$ whose restriction to $U$  
is in $C^\infty(U)$, we wish to show that the restriction
of $u$ to $U$ is in $C^\infty$. Without loss of generality we assume  
that the distribution $u$ has compact support and lies in
$H^{-s_0}(\mathbb R^3)$. For $\delta>0$ we will define a smoothing  
operator $K_\delta$ such that $K_\delta u\in C^\infty$ and
${\rm lim}_{\delta\to 0}K_\delta(\varrho u^+)\sim\varrho u^+$.

\begin{defin} Let $\omega\in\ C_0^\infty(\mathbb R)$, with  
$\omega(0)=1$ and let
$\kappa_\delta(\xi)=\omega(\delta\xi_3)\gamma^+(\xi)$ and
$$
\widehat{K_\delta u}(\xi)=\kappa_\delta(\xi)\hat u(\xi),
$$
where $\gamma^+(\xi)=1$ in a neighborhood of the support of $\hat{u}^+$.
\end{defin}

\begin{lemma} If $\|K_\delta(\varrho u^+)\|_s\le C$ and if  
$\varrho'u^0\in H^s$ then $\varrho u^+\in H^s$.
\end{lemma}
Proof: We have
$$
\|K_\delta(\varrho u^+)-\varrho u^+\|_s\le\|K_\delta((\varrho  
u)^+)-(\varrho u)^+)\|_s+C\|\varrho'u^0\|_s
$$
and
$$
\lim_{\delta\to  
0}(1+|\xi|^2)^{\frac{s}{2}}\omega(\delta\xi_3)\widehat{(\varrho  
u)^+}(\xi))=
(1+|\xi|^2)^{\frac{s}{2}}\widehat{(\varrho u)^+}(\xi)).
$$
Then $(\varrho u)^+\in H^s$ and since $(\widehat{\varrho  
u)^+}-\widehat{\varrho u^+}$ is supported in the elliptic region  
$\mathfrak U^0$ we
have
$$
\|\varrho u^+\|_s\le\|(\varrho u)^+\|_s+C\|\varrho'u^0\|_s ,
$$
thus concluding the proof.
 
\begin{lemma} For $\delta>0${\rm ,} $K_\delta$ is a pseudodifferential  
operator of order $-\infty$ which is of order zero uniformly in  
$\delta$.
$K_\delta$ has the following commutation properties.
\begin{enumerate}
\item[{\rm 1.}]  $[E,K_\delta](I-\Gamma^0)$ is a pseudodifferential operator of  
order $-\infty$ uniformly in $\delta$.
\item[{\rm 2.}] If $R^s$ is a pseudodifferential operator of order $s$ then
$$
[R^s,K_\delta]=\Gamma^0R^{s-1}_\delta+\Psi^{s-1}R^{0}_\delta+R^{- 
\infty}_\delta,
$$
where $R^{s-1}_\delta${\rm ,} $R^{0}_\delta${\rm ,} and $R^{-\infty}_\delta$ are  
pseudodifferential operators of orders $-\infty$ for $\delta>0$ and
of orders $s-1$ and  $0$ uniformly in $\delta$.
\end{enumerate}
\end{lemma}

\Proof  Number 1 follows from the fact that when  
$|\xi|\ge 1$ then $\gamma^0(\xi)=1$ on the support of these symbols. To  
deal
with number 2 we write the principal symbol of $[R^s,K_\delta]$.  
Setting $x_1=x,\ x_2=y$  and $x_3=t$, we have
$$
\sum_j\frac{\partial\kappa_\delta}{\partial\xi_j}\frac{\partial  
r^s}{\partial x_j}=
\delta\omega'(\delta\xi_3)\tilde\gamma^+\frac{\partial r^s}{\partial  
x_3}+
\sum_j\omega(\delta\xi_3)\frac{\partial\tilde\gamma^+}{\partial\xi_j} 
\frac{\partial r^s}{\partial x_j}.
$$
The lemma then follows, since
$$
\delta\omega'(\delta\xi_3)\tilde\gamma^+\frac{\partial r^s}{\partial  
x_3}
=\xi^{s-1}_3\gamma^+\left\{{\tilde\gamma}^+\xi_3^{- 
s}\delta\xi_3\omega'(\delta\xi_3)\frac{\partial r^s}{\partial x_3}\right\},
$$
where ${\tilde\gamma}^+=1$ in a neighborhood of the support of  
$\gamma^+$ and equals zero in a neighborhood of the origin. The
expression in braces is the symbol of an operator of order zero  
uniformly in $\delta$.

\demo{Conclusion of proof of Theorem {\rm C}}  Substituting $K_\delta  
u$ for $u$ in (\ref{final}) we obtain
$$
\|\Psi^{s+\sigma}\varrho K_\delta u^+\|^2\le C(\|\varrho' E_kK_\delta  
u\|_{s+\sigma+k-1}^2+\|\varrho'K_\delta u\|^2_s+
\|K_\delta u\|^2_{-\infty}).
$$
Then we have
\begin{eqnarray*}
\|K_\delta(\varrho u^+)\|_{s+\sigma}^2&\le &C(\|\Psi^{s+\sigma}\varrho  
K_\delta u^+\|^2+\|\varrho'u\|^2_{s+\sigma-1}),\\
 \|[\varrho' E_k,K_\delta] u\|_{s+\sigma+k-1}^2&\le  &
C(\|\varrho''u^0\|^2_{s+\sigma+k-1}+\|u\|^2_{-\infty}),\\
 \|\varrho'K_\delta u\|^2_s&\le& C\|\varrho'u\|^2_s,\\
\noalign{\noindent{\rm and}}
\|K_\delta u\|^2_{-\infty}&\le &C\|u\|^2_{-\infty}.
\end{eqnarray*}
Further
$$
\|\varrho''u^0\|^2_{s+\sigma+k-1}\le  
C(\|\varrho'''E_ku\|^2_{s+\sigma+k-3}+\|u\|^2_{-\infty}).
$$
Therefore, changing notation for the cutoff functions, we get
$$
\|K_\delta(\varrho u^+)\|_{s+\sigma}^2\le C(\|\varrho'  
E_ku\|_{s+\sigma+k-1}^2+\|\varrho'u\|^2_s+\|u\|^2_{-s_0}).
$$
  Therefore, if $u\in H^{-s_0}$, if $u^+\in H^s_{\rm loc}(U)$, and if  
$E_ku\in H^{s+\sigma+k-1}_{\rm loc}(U)$ then $u^+\in H^{s+\sigma}_{\rm loc}(U)$.
It then follows that if $u\in H^{-s_0}$ and if $E_ku\in  
H^{s_1}_{\rm loc}(U)$ then $u^+\in H^{s_1-k+1}_{\rm loc}(U)$. Since, under
the same assumptions, we have $u^0\in H^{s_1+2}_{\rm loc}(U)$ and $u^-\in  
H^{s_1+1}_{\rm loc}(U)$ we conclude that
$u\in H^{s_1-k+1}_{\rm loc}(U)$, thus proving Theorem~C.

\section{Local existence in $L^2$}

The {\it a priori\/} estimates for $E_k$ imply the following local existence result.

\sdemo{Theorem} {\it If $P\in U\subset\mathbb R^3$ with $U$ an open set{\rm ,} then
there exists a neighborhood $U_1\subset\bar U_1\subset U${\rm ,} with $P\in U_1${\rm ,}
 such that if $f\in H^{k-1}_{\rm loc}(U)$
then there exists $u\in L^2(U_1)$ and $E_ku=f$ in $U_1$.}

\Proof   Let $U_1$ be a small neighborhood of $P$. In Lemma 11 set
 $\varrho=1$ in a neighborhood 
of $\bar U_1$ and set $u=v\in C_0^\infty(U_1)$ so that $\varrho u=v$ and 
$[\Psi^{s+\sigma},\Gamma^+]$ is an
operator of order $-\infty$ on $C_0^\infty(U_1)$. Hence we obtain
$$
\|\Psi^{s+\sigma}v^+\|^2\le C(\|E_kv\|^2_{s+\sigma+k}+\|v\|^2_s),
$$
for all $v\in C_0^\infty(U_1)$. Setting $s+\sigma+k=0$ and combinig with the estimates 
for $v^0$ and $v^-$, we 
obtain
$$
\|v\|^2_{-k+1}\le C(\|E_kv\|^2+\|v\|^2_{-k+1-\sigma}).
$$
Then, if the diameter of $U_1$ is sufficiently small, we have
 $$\|v\|^2_{-k+1-\sigma}\le \hbox{small  const. } \|v\|^2_{-k+1}.
$$
Hence
$$
\|v\|_{-k+1}\le \hbox{const. }\|E_kv\|,
$$
for all $v\in C_0^\infty(U_1)$.

  Let $\mathcal W=C_0^\infty(U_1)$ and let $K:\mathcal W\to\mathbb C$ 
be the linear functional defined by
$Kw=(v,f)$ with $w=E_kv$. Then
$$
|Kw|=|(v,f)|\le\|v\|_{-k+1}\|f\|_{k-1}\le C\|w\|.
$$
So $K$ is bounded on $\mathcal W$; hence it can be extended to a bounded linear 
functional on $L^2(U_1)$. Therefore there
exists $u\in L^2(U_1)$ such that $Kw=(w,u)$, that is $(v,f)=(E_kv,u)=(v,E_ku)$. 
Thus $E_ku=f$ in $L^2(U_1)$, which completes the proof.
\vglue-12pt

{
}
 
\centerline{\small (Received August 3, 2003)}

%%%%%%%Theorem info here
%\newtheorem{theorem}{Theorem}[section]
%\newtheorem{proclaim}[theorem]{Proclaim}
 
 %%%%%%%%%%%%%% equation numbering style

%% 
%%%%%%%%%%%%%%%%%%%%
%
 
\def\A {{\mathcal{A}}}
\def\D {{\mathcal{D}}}
\def\R {{\mathbb{R}}}
\def\N {{\mathbb{N}}}
\def\C {{\mathbb{C}}}
\def\Z {{\mathbb{Z}}}
\def\l {\ell}
\def\ml {multline}
\def\multiline {\multline}
\def\lessim {\lesssim}
 
\def\phi{\varphi}
\def\epsilon{\varepsilon}

 \vglue18pt
\begin{center}
{\bf \Large Appendix:}\\
{\bf \Large  Analyticity and loss of
derivatives}  
  \end{center}
  \vglue5pt \centerline{\small By {\scshape Makhlouf Derridj} and {\scshape David S.\
Tartakoff}} 
\renewcommand{\institution}[1]
{\renewcommand{\theinstitutions}{\vskip16pt\baselineskip10pt\begin{quote}
\scriptsize\scshape #1\end{quote}}}

 \institution{5 rue de la Juviniere, 78350 Les Loges en
Josas, France\\
\email{derridj@club-internet.fr}
\\
\vglue-9pt
University
of Illinois at Chicago, Chicago IL\\
\email{dst@uic.edu\\}
  \vglue-30pt \phantom{oh}}

 \shortname{Makhlouf Derridj and David S.\ Tartakoff}

 \shorttitle{Appendix}
\setcounter{section}{0}
\vglue18pt 
\centerline{\bf Abstract} 
\vglue12pt
In \cite{K2005}, J.\ J.\  Kohn
proves $C^\infty$ hypoellipticity for a sum of squares
of complex vector fields which exhibit a large loss of
derivatives. Here, we prove analytic hypoellipticity for
this operator.
  
\section{Introduction and outline}
\renewcommand{\theequation}{\thesection.\arabic{equation}}
\setcounter{equation}{0}
 
In \cite{K2005}, J.\ J.\   Kohn proves hypoellipticity
for the operator 
$$P = LL^*+ (\overline{z}^kL)^*(\overline{z}^kL),
\qquad L = {\partial  \over \partial z} + i\overline{z}
{\partial \over \partial t},$$
for which there is a large loss of
derivatives --- indeed in the {\it a priori} estimate
one  bounds only the Sobolev norm of order $-(k-1)/2$,
and thus there is a loss of $k-1$ derivatives: $Pu\in
H^s_{\rm loc}\implies u\in H^{s-(k-1)}_{\rm loc}$.

We show in this note that solutions of
$Pu=f$ with $f$ real analytic are themselves real
analytic in any open set where $f$ is. In so doing we
 use an {\it a priori} estimate which follows easily
from that established by Kohn for this operator,
namely for test functions $v$ of small support near
the origin: 
\begin{equation}\label{apeDT}\|\overline{L}v\|_0^2 +
\|\overline{z}^k{L}v\|_0^2 + \|v\|^2_{-{k-1\over 2}}
\lesssim |(Pv, v)_{L^2}|.\end{equation}
In fact, in \cite{T3} (see also \cite{BDKT}),
we give a rapid and direct derivation
of (\ref{apeDT}) for this operator and similar
estimates for more degenerate operators.

The first two
terms on the left of this estimate exhibit maximal
control in
$\overline{L}$ and $\overline{z}^kL$, but only these
complex directions. Hence in obtaining recursive
bounds for derivatives it is essential to keep one of
these vector fields available for as long as possible.
For this, we will construct a carefully balanced
localization of high powers of $T= -2i
\partial/\partial t$  and use the estimate repeatedly,
reducing the order of powers of $T$ but accumulating
derivatives on the localizing functions. These
Ehrenpreis type localizing functions work `as
if analytic' up to a prescribed order, with
all constants independent of that order, as in 
\cite{DT1978}, \cite{DT1980}, but eventually the good
derivatives ($\overline{L}$ or
$\overline{z}^kL$) are lost and we must use the third
term on the left of the estimate, absorb the
loss of $\frac{k-1}{2}$ derivatives, introduce a new
localizing function of larger support and start the
whole process again, but with only a (fixed) fraction
of the original power of $T$.

\section{Observations and simplifications}

Our first observation is that we know the analyticity
of the solution for $z$ different from $0$ from the
earlier work of the second author \cite{DT1978},
\cite{DT1980} and Tr\`eves \cite{Tr1978}. Thus, modulo
brackets with localizing functions whose derivatives
are supported in the known analytic hypoelliptic
region, we take all localizing functions independent of
$z$.

Our second observation is that it suffices to bound
derivatives measured in terms of high powers of the
vector fields $L$ and $\overline{L}$ in $L^2$ norm, by
standard arguments, and indeed estimating high powers
of ${L}$ can be reduced to bounding high
powers of 
$\overline L$ and powers of $T$ of half the order, by
repeated integration by parts.  Thus our overall scheme
will be to start with high powers (order $2p$) of $L$ or
$\overline{L},$ use integration by parts and
the {\it a priori} estimate repeatedly to reduce to
treating $T^pu$ in a slightly larger set.  

And to do this, we introduce a new special
localization of
$T^p$ adapted to this problem.

\section{The localization of high powers of $T$}

The new localization of $T^p$ may be written in the
form: 
$$(T^{p_1,p_2})_\phi = \sum_{a\leq p_1 \atop b\leq
p_2}{L^a\circ 
z^a\circ 
T^{p_1-a}\circ
\phi^{(a+b)}\circ T^{p_2-b}\circ
\overline{z}^b\circ \overline{L}^b
\over a!b!}.$$
Here by $\phi^{(r)}$ we mean $(-i\partial/\partial
t)^r\phi(t)$ since near $z=0$ we have seen that we may
take the localizing function independent of $z.$ Note
that the leading term (with $a+b=0)$ is merely
$T^{p_1}\phi T^{p_2}$ which equals $T^{p_1+p_2}$
on the initial open set $\Omega_0$ where $\phi \equiv
1$.

We have the commutation relations: 
$$[L, (T^{p_1,p_2})_\phi] \equiv L\circ
(T^{p_1-1,p_2})_{\phi'},$$
$$ [\overline{L},
(T^{p_1,p_2})_\phi] \equiv 
(T^{p_1,p_2-1})_{\phi'}\circ \overline{L},
$$
$$[(T^{p_1,p_2})_\phi,z] = (T^{p_1-1,p_2})_{\phi'}\circ
z,
$$
and
$$
[(T^{p_1,p_2})_\phi,\overline{z}] =
\overline{z}\circ (T^{p_1,p_2-1})_{\phi'},
$$
\vskip.1in\noindent
where the $\equiv$ denotes modulo
$C^{p_1-p_1'+p_2-p_2'}$ terms of the form
\begin{equation}\label{pureL}{L^{p_1-p_1'}\circ
z^{p_1-p_1'}\circ 
T^{p_1'}\circ\phi^{(p_1-p_1'+p_2-p_2'+1)}\circ T^{p_2'}
\circ\overline{z}^{p_2-p_2'}\circ
\overline{L}^{p_2-p_2'}\over (p_1-p_1')!(p_2-p_2')!}
\end{equation}
with either $p_1'=0$ or
$p_2'=0$, i.e., terms where all free $T$
derivatives have been eliminated on one side of
$\phi$ or the other. Thus if we start with
$p_1=p_2={p/ 2}$, and iteratively apply these
commutation relations, the number of $T$ derivatives not
necessarily applied to
$\phi$ is eventually at most ${p/ 2}$. 

\section{The recursion}

We insert first $v=(T^{{p\over 2},{p\over 2}})_\phi u$
in the {\it a priori} inequality, then bring
$(T^{{p\over 2},{p\over 2}})_\phi$ to the left of
$P=-L\overline{L}-
 \overline{L}z^k\overline{z}^kL$ since $Pu$ is known
and analytic. We have, omitting for now the
`subelliptic' term,
$$\|\overline{L}(T^{{p\over 2},{p\over 2}})_\phi
u\|_0^2 +
\|\overline{z}^k{L}(T^{{p\over 2},{p\over 2}})_\phi
u\|_0^2 \lesssim |(P(T^{{p\over 2},{p\over 2}})_\phi u,
(T^{{p\over 2},{p\over 2}})_\phi u)_{L^2}|$$
$$\lesssim
|((T^{{p\over 2},{p\over 2}})_\phi Pu, (T^{{p\over 2},{p\over 2}})_\phi
u)_{L^2}| + |([P,(T^{{p\over 2},{p\over 2}})_\phi] u,
(T^{{p\over 2},{p\over 2}})_\phi u)_{L^2}| 
$$
and, by the above bracket relations,  
\begin{eqnarray*}
&&([P,(T^{{p\over 2},{p\over 2}})_\phi] u,
(T^{{p\over 2},{p\over 2}})_\phi u)\\
&&\qquad = -([L\overline{L},(T^{{p\over 2},{p\over 2}})_\phi]
u, (T^{{p\over 2},{p\over 2}})_\phi u) - ([
\overline L z^k\overline{z}^k{L},(T^{{p\over 2},{p\over
2}})_\phi] u, (T^{{p\over 2},{p\over 2}})_\phi u)\\
&&\qquad   \equiv -(L(T^{{{p\over 2},{p\over 2}}-1})_{\phi'}
\overline{L}u, (T^{{p\over 2},{p\over 2}})_\phi u)
-(L (T^{{p\over 2}-1,{p\over 2}})_{\phi'}
\overline{L}u, (T^{{p\over 2},{p\over 2}})_\phi u)\\
&&\qquad \quad-
((T^{{p\over 2}-1,{p\over
2}})_{\phi'}\overline{L} z^k\overline{z}^k {L}u, 
(T^{{p\over 2},{p\over 2}})_\phi u) 
\\
&&\qquad \quad-\sum_{k'=1}^k 
(\overline L{z}^{k'}(T^{{p\over 2},{p\over
2}-1})_{\phi'}
{z}^{k-k'}\overline z^k{L}u, 
(T^{{p\over 2},{p\over 2}})_\phi u) \\
&&\qquad \quad-\sum_{k'=0}^{k-1}(\overline L{z}^k
\overline z^{k'} (T^{{p\over 2}-1,{p\over 2}})_{\phi'}
\overline z^{k-k'}{L}u, 
(T^{{p\over 2},{p\over 2}})_\phi u)\\
&&\qquad \quad- (\overline L{z}^k
\overline z^k{L}(T^{{p\over 2},{p\over
2}-1})_{\phi'}u,  (T^{{p\over 2},{p\over
2}})_\phi u),
\end{eqnarray*}
with the same meaning for $\equiv$ as above. In every
term, no powers of $z$ or $\overline{z}$ have been
lost, though some may need to be brought to the left
of the $(T^{q_1,q_2})_{\tilde{\phi}}$ 
with again no loss
of powers of $z$ or $\overline{z}$ and a further
reduction in order, every bracket reduces the order of 
the sum of the two indices $p_1$ and
$p_2$ by one (here we started with $p_1=p_2=p/2$),
picks up one derivative on $\phi,$ and leave the vector
fields over which we have maximal control in the
estimate intact and in the correct order. Thus we may
bring either
$\overline{L}z^k$ or $L$ to the right as
 $\overline{z}^kL$ or
$\overline{L},$ and use a weighted Schwarz inequality on
the result to take maximal advantage of the {\it a
priori} inequality. Iterations of all of this
continue until there remain at most $p/2$
free $T$ derivatives (i.e., the $T$ derivatives on
at least one side of $\phi$ are all `corrected' by
good vector fields) and perhaps as many as
${p/2} \;L$ or $\overline{L}$ derivatives, and we
may continue further until, at worst, these remaining
$L$ or $\overline{L}$ derivatives bracket two at a
time to produce more $T$'s, with corresponding
combinatorial factors. After all of this, there will be
at most
$T^{3p\over 4}$ remaining, and a factor of
${p\over 2}!!\sim {p\over 4}!$ 

It is here that the final term on the left of the
{\it a priori} inequality is used, in order to bring
the localizing function out of the norm after
creating another balanced localization of $T^{3p/4}$
with a new localizing
function of Ehrenpreis type with slightly larger
support, geared, roughly, to
$3p/4$ instead of to
$p$. 

Recall that such such localizing functions
$\psi$ may be constructed for any $N$ and satisfy 
$$\left|\psi^{(r)}\right| \leq \left({C\over
e}\right)^{r+1}N^r,
\quad r\leq 2N$$
where $C$ is independent of $N$ and 
$e={\rm dist}(\{\psi\equiv 1\}, ({\rm supp}\,\psi)^c)$.

\section{Conclusion of the proof}

Finally, this entire process, which reduced the
order from $p$ to at most  $3p/4,$ (or more precisely to at most $3p/4 +
(k-1)/2$), is repeated, over and over, each time
essentially reducing the order by a factor of $3/4.$
After at most  $\log_{4/3} p$ such iterations
we are reduced to a bounded number of derivatives,
and, as in \cite{DT1978} and \cite{DT1980}, all of these
nested open sets may be chosen to fit in the one open
set $\Omega'$ where $Pu$ is known to be analytic, and
all constants chosen independent of $p$ (but
depending on $Pu$). The fact that in those
works one full iteration reduced the order by half
played no essential role --- a factor of
$3/4$ works just as well. 

To be precise, the sequence of open sets,
$\{\Omega_j\}$, each compactly contained in the next,
with
$\Omega_{\log_{4/3}p}=\Omega'$, have separations
$d_j={\rm dist}(\Omega_j, \Omega_{j+1}^c),$ with $\sum
d_j={\rm dist}(\Omega_0,{\Omega'}^c)=d$, which need to be
picked carefully. The localizing functions
$\{\phi_j\}$ with $\phi_j\in C_0^\infty(\Omega_{j+1}) 
\equiv 1$ on
$\Omega_j$
satisfy 
\begin{equation}\label{j}\left|\phi_j^{(r)}\right|\leq
(C/d_j)^{r+1}((3/4)^jp)^r, \qquad r\leq 2(3/4)^jp.
\end{equation}
We shall take the $d_j={1\over {(j+1)^2}}/d\sum {1\over
{(j+1)^2}}$. 
\smallbreak 
Now at most $(3/4)^jp$ derivatives will fall on
$\phi_j$, and most of the effect of the derivatives
will be balanced by corresponding factorials in the
denominator, as in (\ref{pureL}), roughly the powers of
$(3/4)^jp$ in (\ref{j}) in view of Stirling's formula.
In addition, as noted immediately before the last
paragraph in Section 4, there will be factorials
corresponding to the diminution of powers of $T$. What
will {\it not} be balanced are the powers of
$d_j^{-1}$, but the product of these factors will
contribute 
$$\Pi_{j=1}^{\log_{4/3}p}
\left(j^2\right)^{(3/4)^jp}= 
\left(\Pi_{j=1}^{\log_{4/3}p}
j^{(3/4)^j}\right)^{2p}= C^p,$$
which, together with the factorials just mentioned,
proves the analyticity of the solution in
$\Omega_0$.  

 \vglue-6pt

\end{document}